\documentclass[11pt]{amsart}

\usepackage[english]{babel}
\usepackage{amsmath}
\usepackage{amsfonts}
\usepackage{amssymb}
\usepackage{amsopn}
\usepackage{amsthm}
\usepackage{enumerate}
\usepackage{dsfont}
\usepackage{mathrsfs}
\usepackage{extarrows}
\usepackage{wasysym}
\usepackage{stmaryrd}
\usepackage{todonotes}
\usepackage{algpseudocode}
\usepackage[shortlabels]{enumitem}
\usepackage{xspace}
\usepackage[arrow,matrix,curve]{xy} 
\usepackage{color, tikz}
\usepackage{listings}
\usepackage{hyperref}
\usepackage[nameinlink]{cleveref}

\usepackage{tablefootnote}
\usepackage{multirow}
\usepackage{xcolor,colortbl}

\usepackage{siunitx}
\usepackage{caption}
\usepackage{subcaption}
\usepackage{float}
\usepackage{bbm}
\usepackage[foot]{amsaddr}

\newcounter{FNC}[page]
\def\newfootnote#1{{\addtocounter{FNC}{2}$^\fnsymbol{FNC}$%
     \let\thefootnote\relax\footnotetext{$^\fnsymbol{FNC}$#1}}}

\crefname{remark}{remark}{remarks}
\crefformat{remark}{#2Remark~#1#3}

\crefname{section}{section}{Sections}
\crefformat{section}{#2Section~#1#3}

\crefname{subsection}{Subsection}{Subsections}
\crefformat{subsection}{#2Section~#1#3}

\crefname{appendix}{Appendix}{Appendices}
\crefformat{appendix}{#2Appendix~#1#3}

\crefname{equation}{equation}{equations}
\crefformat{equation}{#2(#1)#3}

\crefname{enumi}{enumi}{enumis}
\crefformat{enumi}{#2#1#3}

\crefname{figure}{Figure}{Figures}
\crefformat{figure}{#2Figure~#1#3}

\crefname{table}{Table}{Tables}
\crefformat{table}{#2Table~#1#3}

\numberwithin{equation}{section}
\numberwithin{environment}{section}

\makeatletter
\def\subsubsection{\@startsection{subsubsection}{3}%
  \z@{.5\linespacing\@plus.7\linespacing}{-.5em}%
  {\normalfont\bfseries}}
\makeatother

\def\Id{\mathbbm{1}}

\def\O{\mathcal{O}}
\def\N{\mathcal{N}}

\newcommand{\norm}[1]{\left\lVert#1\right\rVert}

\DeclareMathOperator{\tr}{tr}
\DeclareMathOperator{\sgn}{sign}
\DeclareMathOperator{\diag}{diag}

\newcommand{\Vector}[1]{\ensuremath{\boldsymbol{#1}}}

\definecolor{NiceBlue}{rgb}{0.2,0.2,0.75}
\newif\ifstartedinmathmode
\newcommand{\struc}[1]{{\relax\ifmmode\startedinmathmodetrue\else\startedinmathmodefalse\fi\color{NiceBlue}{\ifstartedinmathmode #1 \else\textit{#1}\fi}}}

\newcommand{\tHU}{\ensuremath{t_{\text{HU}}}}
\newcommand{\tSOS}{\ensuremath{t_{\text{SOS}}}}
\newcommand{\tHUQuantum}{\ensuremath{t_{\text{HU}}^{\text{Quantum}}}}
\newcommand{\zHU}{\ensuremath{Z_{\text{HU}}}}
\newcommand{\zSOS}{\ensuremath{Z_{\text{SOS}}}}
\newcommand{\DeltaHU}{\ensuremath{\Delta_{\text{HU}}^{\text{rel}}}}
\newcommand{\DeltaSOS}{\ensuremath{\Delta_{\text{SOS}}^{\text{rel}}}}
\newcommand{\zQUBO}{\ensuremath{Z_{\text{QUBO}}}}
\newcommand{\DeltaQUBO}{\ensuremath{\Delta_{\text{QUBO}}^{\text{rel}}}}
\newcommand{\DeltaHUabs}{\ensuremath{\Delta_{\text{HU}}^{\text{abs}}}}
\newcommand{\DeltaSOSabs}{\ensuremath{\Delta_{\text{SOS}}^{\text{abs}}}}
\newcommand{\DeltaQUBOabs}{\ensuremath{\Delta_{\text{QUBO}}^{\text{abs}}}}
\newcommand{\zHUupper}{\ensuremath{\widetilde{Z}_{\text{HU}}}}
\newcommand{\DeltaHUupper}{\ensuremath{\widetilde{\Delta}_{\text{HU}}^{\text{rel}}}}
\newcommand{\DeltaHUabsupper}{\ensuremath{\widetilde{\Delta}_{\text{HU}}^{\text{abs}}}}

\newcommand{\zVRP}{\ensuremath{Z_{\text{OVRP}}^*}}
\newcommand{\Zopt}{\ensuremath{Z^*}}

\title[Benchmarking SDP relaxations of QUBO formulations]{Benchmarking of quantum and classical SDP relaxations for QUBO formulations of real-world logistics problems}
\author[]{Birte Ostermann$^{1}$, Taylor Garnowski$^{\dagger,3}$, Fabian Henze$^{\dagger,2}$, Vaibhavnath Jha$^{\dagger,4}$, Asra Dia$^{1}$, Frederik Fiand$^{4}$, David Gross$^{2}$, Wendelin Gross$^{3}$, Julian Nowak$^{3}$, \and Timo de Wolff$^{*,1}$}

\address{$^{\dagger}$ These authors contributed equally to this work.}
\address{$^*$ Corresponding author: \textit{t.de-wolff@tu-braunschweig.de}}

\address{$^1$ Institute for Analysis and Algebra, TU Braunschweig, Germany}
\address{$^2$ Institute for Theoretical Physics, University of Cologne, Germany}
\address{$^3$ 4flow SE, Hallerstrasse 1,  10587 Berlin, Germany}
\address{$^4$ GAMS Software GmbH, PO Box 4059, 50216 Frechen, Germany}



\begin{document}
\date{}

\begin{abstract}
Quadratic unconstrained binary optimization problems (QUBOs) are intensively discussed in the realm of quantum computing and polynomial optimization.
We provide a vast experimental study of semidefinite programming (SDP) relaxations of QUBOs using sums of squares methods and on Hamiltonian Updates.
We test on QUBO reformulations of industry-based instances of the (open) vehicle routing problem and the (affinity-based) slotting problem -- two common combinatorial optimization problems in logistics.
Beyond comparing the performance of various methods and software, our results reaffirm that optimizing over non-generic, real-world instances provides additional challenges. In consequence, this study underscores recent developments towards structure exploitation and specialized solver development for the used methods and simultaneously shows that further research is necessary in this direction both on the classical and the quantum side.
\end{abstract}

\maketitle

\section{Introduction}
\textbf{Problem Setting and Approach.} In the last 25 years there has been plenty of development on solving nonlinear optimization problems and especially \struc{(constrained) polynomial optimization problems (CPOPs)}
\begin{equation*}
		\begin{aligned}
			&\text{min}							
			& &	f(\Vector{x})\\
			&\text{subject to}	&	&	h_1(\Vector{x}) \geq 0,\ldots, h_s(\Vector{x}) \geq 0 \\
			&	&	&	\Vector{x} \in \mathbb{R}^n,
		\end{aligned}
	\end{equation*}
where the objective function $f(\Vector{x})\in \mathbb{R}[\Vector{x}] = \mathbb{R}[x_1,\ldots,x_n]$ and the constraints $h_1(\Vector{x}),\dots,h_s(\Vector{x})\in\mathbb{R}[\Vector{x}]$ are polynomial.
A prominent subclass of CPOPs consists of \struc{quadratic unconstrained binary optimization (QUBO)} problems, which are of the form
\begin{align*}
		\min_{\Vector{x}\in \{0,1\}^{n}}\Vector{x}^TQ\Vector{x},
\end{align*}
for a given real symmetric $n \times n$ matrix $Q\in \mathbb{R}^{n\times n}$.

QUBOs received a lot of attention not only because one can represent several real-world applications by them, but also due to their connection to approaches based on quantum computing.
In fact, there exists a wide range of algorithms based on different classical and quantum (inspired) ans\"atze for solving QUBOs.
Prominently, but not exclusively, these include 
\begin{enumerate}
	\item interpreting the QUBO as a CPOP, which is tackled via \struc{sums of squares (SOS)}, or, on the dual side, \struc{moment methods}. This approach leads to a hierarchy or lower bounds, known as \struc{Lasserre's hierarchy} \cite{lasserre, parrilo2000structured}.
	Both the primal and the dual approach induce a \struc{semidefinite program (SDP)}, which can classically be solved, e.g.~using \struc{interior point methods}.
    \item similarly, form the quantum community perspective, interpreting the QUBO as in the Goemans-Williamson \cite{Goeman_Williamson} approximation algorithm for the \struc{maximum cut} \struc{(MaxCut)} problem.
    This also results in an SDP (which coincides with the SDP in the first order of the Lasserre hierarchy), that is solved using \struc{Hamiltonian Updates (HU)} \cite{Brandao2022fasterquantum,henze2025}, providing a classical and a quantum version of the algorithm.
	\item solving the QUBO directly without an intermediate relaxation via classical, specialized QUBO solvers, e.g.~\cite{QuBowl, McSparse}, or \struc{quantum annealing} \cite{DWaveLeap}.
\end{enumerate}
All of these methods have been implemented and there exist various experimental evaluations of these methods.
In this article, we focus in particular (but not exclusively) on a comparison of the 
methods listed in the first two points, i.e.~the SOS based methods, referred to as \struc{SOS-SDP}, and HU.

\medskip

In this interdisciplinary study we provide a vast experimental comparison of a broad range of methods and solvers for solving SDP relaxations of QUBOs formulations of two common problems in logistics.
We choose these sets of problems, as applications in science and engineering regularly lead to highly non-generic optimization problems, which are known to behave differently than randomly generated data.
Specifically, we investigate the \struc{(open) vehicle routing problem (OVRP)} concerning the distribution of goods carried out with a fleet of transporters, and the \struc{(affinity-based) slotting problem (ASP)} concerning the optimal storage of goods in warehouses.
The problems in our study are represented by \struc{real-world data} -- in opposition to using generated, random instances -- which is provided by the industry partner 4flow SE.
We chose the OVRP and ASP instances not only because  they can be reformulated as QUBOs, their practical relevance, and the availability of real-world data, but also because they are natural candidates for benchmarking: 
The original mathematical representation of the problems are \struc{integer programs (IP)} and \struc{integer quadratic programs (IQP)}, respectively, so that IP/IQP solvers can still effectively solve the specific instances and thus provide benchmarks.
This allows us to not only compare the different approaches among each other, but also make qualitative statements about their effective quality on these problem sets.

In terms of the SDP methods, we focus in particular on different software and solvers for SOS-SDP based hierarchies, and on HU.
We highlight that we do not intend to propose a specific approach, but give a qualitative comparison and highlight abilities and challenges for the different methods.
All results, logs, and a full documentation of the experiments are available via

\smallskip
\url{https://moto.math.nat.tu-bs.de/appliedalgebra_public/sdp_qubo_benchmarking}.

\medskip

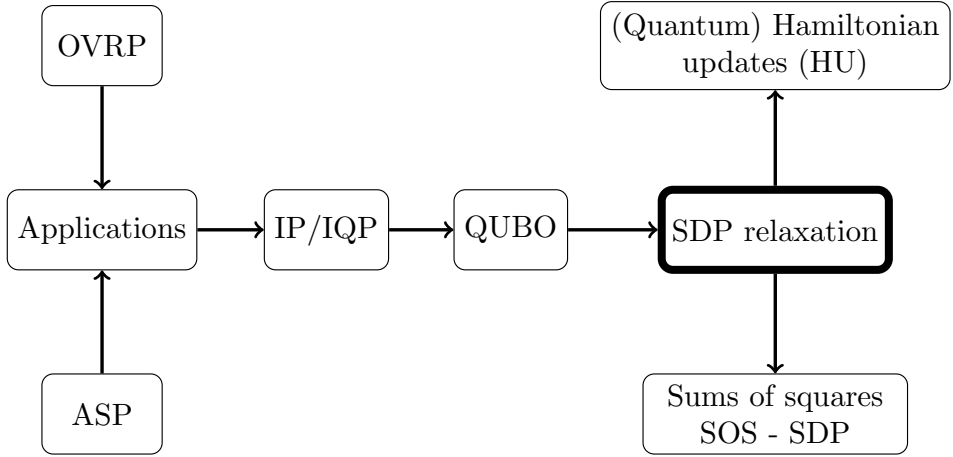
\begin{figure}[h] 
			\centering
			\resizebox{0.8\linewidth}{!}{
				\begin{tikzpicture}[node distance=2.3cm]
					\node (app) [rectangle, rounded corners,  minimum height=1cm,text centered, draw=black, fill=white] {Applications};
					\node (VRP) [rectangle, rounded corners, minimum width=1.5cm, minimum height=1cm,text centered, draw=black, fill=white, above of=app, align=center] {OVRP};
					\node (ASP) [rectangle, rounded corners, minimum width=1.5cm, minimum height=1cm,text centered, draw=black, fill=white, below of=app, align=center] {ASP};
					\node (MIP) [ rectangle, rounded corners,  minimum height=1cm,text centered, draw=black, fill=white, right of=app, align=center, xshift = 0.5cm, text centered] {IP/IQP};
			        \node (QUBO) [ rectangle, rounded corners, minimum height=1cm,text centered,draw = black,  fill=white, right of=MIP, align=left] {QUBO};
					\node (relax) [ line width = 1mm, rectangle, rounded corners, minimum height=1cm, text centered, draw=black, fill=white, right of=QUBO, align=left, xshift=1cm] {SDP relaxation};
					\node (qSDP) [ rectangle, rounded corners, minimum width = 3.3cm, minimum height=1cm,text centered, draw=black, fill=white, above of=relax, align=center] {(Quantum) Hamiltonian\\ updates (HU)};
					\node (SOS) [ rectangle, rounded corners,  minimum width = 3.3cm, minimum height=1cm, text centered, draw=black, fill=white, below of=relax, align=center] {Sums of squares\\SOS - SDP};
					\draw [black, ->, very thick, align=left] (VRP) -- node[anchor=south] {} (app);
					\draw [black, ->, very thick, align=left] (ASP) -- node[anchor=south] {} (app);
					\draw [black, ->, very thick, align=left] (app) -- node[anchor=south] {} (MIP);
					\draw [black, ->, very thick, align=left] (MIP) -- node[anchor=south] {} (QUBO);
					\draw [black, ->, very thick, align=left] (QUBO) -- node[anchor=south] {} (relax);
					\draw [black, ->, very thick, align=left] (relax) -- node[anchor=south] {} (qSDP);
					\draw [black, ->, very thick, align=left] (relax) -- node[anchor=south] {} (SOS);
				\end{tikzpicture}
			}
		\caption{Visualization of the pipeline for this work: We compare different methods and software for semidefinite programming (SDP) relaxations of quadratic unconstrained binary optimization (QUBO) formulations that emerge from real-life applications. In particular, we consider the open vehicle routing problem (OVRP) and the affinity-based slotting problem (ASP). The applications are initially formulated as integer program (IP) or integer quadratic program (IQP) with either a linear objective function (for OVRP) or a quadratic objective function (for ASP), respectively, that can be solved to optimality.
        }
        \label{fig:pipeline}
	\end{figure}

\textbf{Results.} Our main results are:
\begin{enumerate}
	\item Due to the IP/IQP structure of the original problem instances, our approach requires a two-step reformulation: First, we reformulate the instances as a QUBO, which we then relax as an SDP, see \cref{fig:pipeline}.
	Hence, we expected that the chosen problems are challenging for all chosen methods.
	The experiments confirm this expectation, showing that all methods and solvers that are not directly working on the IP/IQP side have significant problems with solving the corresponding SDPs; see \cref{sec:results} and in particular \cref{tab:VRPresults} and \cref{fig:ASP_running_times}.
	\item Despite the general challenges, we see significant differences in terms of quality of bounds and running times of the different methods.
	In particular, we see that
	\begin{itemize}
		\item on the SDP-SOS side there is a broad range of qualitative differences.
		In several cases, we observed that the SDP solvers had both memory and numerical problems.
		In terms of quality of the results, the software \textsc{TSSOS} \cite{TSSOS}, which exploits sparsity, joint with the \textsc{Mosek} \cite{mosek} SDP solver obtained the best results.
		\item HU does not outperform the best SOS-SDP results in terms of quality.
		However, the running time of HU is mainly based on the chosen precision, which allows both to tackle larger size problems, and to compute very rough bounds quickly.
		Moreover, in contrast to SOS-SDP, it always provides solution recovery.
	\end{itemize}
    \item  Recently, SDP solvers as well as software for preprocessing were developed that specifically exploit the problem's properties, e.g.~sparsity, low rank or matrix structures. 
    In summary, we confirm the potential for improvement given by the application of such  specialized software.
    Additionally, we see that methods can perform differently or worse on real-world data, showing the importance of including non-randomly generated instances in software developments and benchmark studies.
    \item
    This finding is mirrored on the quantum side.
    Quantum algorithms research is facing the issue that benchmarks are difficult to perform, due to the absence of larger-scale quantum hardware.
    By necessity then, results on the performance of quantum algorithms typically rely on strong assumptions and often pertain only to asymptotic behavior of generic approaches.
    What we see here is that the assessment of quantum-algorithmic performance for real-world problems 
    is not adequately captured by these high-level arguments.
    Practitioners who want to gauge the impact of quantum computing for their use cases will find that the running time of quantum solutions for their specific instances may fall widely short of promises based on asymptotics.
\end{enumerate}

\medskip

\textbf{Organization.} This article is organized as follows: After brief preliminaries in \cref{sec:preliminaries}, we introduce background on our two problem classes, OVRP and ASP, in \cref{sec:Vehicle_Routing_Problems} and \cref{sec:slotting_problems}, respectively, and explain their reformulation from IP/IQP to QUBO structure in \cref{sec:QUBO_reformulation}.
In \cref{sec:SOSSDP} we provide background on the SOS-SDP relaxation, and in \cref{sec:HU} on the Hamiltonian Updates algorithm both from a classical and a quantum perspective.
In the main \cref{sec:experiments} we present the experimental results: We recall the experimental setup in \cref{sec:ex_setup}, provide benchmarks from IP/IQP and QUBO formulation in \cref{sec:optimal_mip_qubo_results}, and the results from the SDP relaxations in \cref{sec:results}, followed by a discussion in \cref{sec:discussion}.
We conclude the paper in \cref{sec:conclusion}.
Several further technical details are specified in \cref{appendix:details_problem_instances,sec:nonnegativity_certificates,appendix:tables,subsection:binary_conversion,app:tech_details}.

\medskip

\textbf{Related Work.} There exists a large range of articles related to this study.
First, regarding the SOS-SDP relaxation there are various general introductions e.g.~\cite{Laurent:Survey, Lasserre:IntroductionPolynomialandSemiAlgebraicOptimization, Theobald:Book:RealAlgGeom}, and in particular several works on using Lasserre's hierarchy on binary optimization problems and their convergence guarantees \cite{LauConjecture, Fawzi2016, sakaue}.

In the context of CPOPs and QUBOs various colleagues made significant effort in providing experimental studies and comparisons. 
There exists a collection of benchmark problem sets and corresponding studies, e.g.~\cite{wiegele2007biq,qplib, mittelmann}.

Related to combinatorial optimization problems discussed here, there are works benchmarking MaxCut algorithms, including the approximation algorithm by Goemans and Williamson, on instances, among other, from real-world graphs from the network library \cite{maxcutbenchmark}.
Moreover, SDP relaxations were applied to various other applications resulting in benchmarking articles. 
A famous recent example is the application of sparsity variants of the SOS based hierarchies to the Alternating-Current Optimal Power Flow Problem using different solvers, e.g.~\textsc{TSSOS} \cite{TSSOSDok,OptPowerflow}, and \textsc{sparsePOP} \cite{OPF}.
Another example is the overview \cite{Ahmadi2016} displaying SOS-SDP (and SDSOS) relaxations for polynomial optimization problems from Operations Research and Transportation Engineering.

In this article we use a wide range of solvers for our experiments. 
There are several works benchmarking these solvers on random instances:
\cite{habibi:hal-04076510} developed the \textsc{Loraine} solver and tested it on randomly generated MaxCut and QUBO instances with $n=50$ variables, comparing it to the standard, commercial SDP solver \textsc{Mosek} \cite{mosek_solver,mosek}.
Similarly, \textsc{ManiSDP} \cite{wang2024solvinglowranksemidefiniteprograms} solved second order Lasserre relaxations on randomly generated QUBOs with up to $n=120$ variables, likewise comparing the results to \textsc{Mosek}.
The solver \textsc{COSMO} \cite{cosmo} is benchmarked on several randomly generated SDP instances.

Hamiltonian Updates was developed in \cite{Brandao2022fasterquantum}. It was improved and benchmarked on randomly generated instances in \cite{henze2025}. The quantum running time analysis uses a similar gate count approach as \cite{simplex2023,Dalzell2023,Babbush_2021, Campbell_2019}.

Furthermore, beyond SDP based methods, \cite{juengerDwave} tested QUBO problems with D-Wave's quantum annealing procedure against branch-and-bound and SDP methods, reporting better results than the classical methods on these particular methods, and \cite{HTX2023} compared quantum annealing and quantum inspired annealers on quadratic assignment problems.

\section{Preliminaries}
\label{sec:preliminaries}
We give an overview on the basic notation and definitions that are essential in the upcoming sections.
The method specific definitions each appear in the beginning of the respective \cref{sec:industryproblems,sec:SOSSDP,sec:HU}.
We refer to the set of \struc{natural}, \struc{integer} and \struc{real} numbers by \struc{$\mathbb{N}$}, \struc{$\mathbb{Z}$} and \struc{$\mathbb{R}$} respectively, and thus denote $\struc{\mathbb{Z}_{\geq 0}} = \{x \in \mathbb{Z}\;|\; x\geq 0\}$ and $\struc{\mathbb{R}_{\geq 0}}= \{x \in \mathbb{R}\;|\; x\geq 0\}$.
We use bold letters to denote vectors, e.g.~$\Vector{x} = [x_1\;\hdots\; x_n]^T \in \mathbb{R}^n$.
With the letter $T$ we indicate the transpose of a vector or matrix. 
In this article, we specifically consider \struc{real} \struc{symmetric} matrices $\struc{A~=~(a_{ij})_{1\leq i,j\leq n}\in \mathbb{R}^{n\times n}}$ with $\struc{A^T = A}$.
A symmetric matrix $A \in \mathbb{R}^{n\times n}$ is \struc{positive semidefinite} if and only if $\Vector{x}^TA\mathbf{x}\geq 0$ for all $\Vector{x}\in \mathbb{R}^n$, and we indicate this via the notation $A\;\struc{\succcurlyeq}\;0$. 
For two matrices $A,B\in \mathbb{R}^{n\times n}$ we define their \struc{inner product} $\struc{\langle A, B\rangle} = \text{trace}(A^TB) = \sum_{i,j}^{n}a_{ij}b_{ij}$.

Let $r\in \mathbb{N}$ and $C, A_1,\dots ,A_r \in \mathbb{R}^{n\times n}$ be symmetric matrices and $\Vector{b}\in \mathbb{R}^r$. 
A \struc{semidefinite optimization problem (SDP)} is of the following form

\begin{equation*}
\begin{aligned}
			&\min_{X\in \mathbb{R}^{n\times n}}						
			& &	\langle C, X\rangle \\
			&\text{subject to}	& &	\langle A_i, X\rangle = b_i\; \text{for all }i=1,\dots,r \\
			&	&	&	X \text{ symmetric},\;X\succcurlyeq 0.
		\end{aligned}\label{eq:SDPgeneral}
\end{equation*}
In this article, we consider SDPs that emerge as relaxations of \struc{quadratic unconstrained binary optimization problems (QUBOs)}. 
Let $Q\in \mathbb{R}^{n\times n}$ be a symmetric matrix.
Then a QUBO is of the form
\begin{equation}
\begin{aligned}
			&\min_{\Vector{x}\in \{0,1\}^{n}}						
			& &	\Vector{x}^TQ\Vector{x} 
		\end{aligned}\label{eq:QUBO}
\end{equation}
In \eqref{eq:QUBO}, we define the QUBO via optimizing over the set $\{0,1\}^n$.
Alternatively, one can define a QUBO over $\{1,-1\}^n$.
This formulation is often called \struc{Ising model} in the physics context and we use it mostly in \cref{sec:HU}.
We give the conversion between both formulations in \cref{subsection:binary_conversion}.
    
\section{Industry Problems and QUBO Reformulations}
\label{sec:industryproblems}
In this section, we introduce two business-relevant combinatorial optimization problems.
Our industry partner, 4flow SE, provided data stemming from their work with customers on real-world problems. 
In what follows, we refer to this data by \struc{real-world data}. 
In particular, we consider a \struc{vehicle routing problem (VRP)} and an \struc{affinity-based slotting problem (ASP)}. 

\subsection{Vehicle Routing Problems}
\label{sec:Vehicle_Routing_Problems}
The \struc{vehicle routing problem (VRP)} is an extension of the well-studied and NP-hard \struc{traveling salesman problem (TSP)}.
State-of-the-art heuristics, such as \cite{helsgaun2000effective}, and exact solvers like Concorde \cite{concorde} can effectively solve TSP instances involving thousands of nodes. 
However, these approaches are highly specialized for the TSP and are generally not applicable in industrial settings. 
Generic TSPs rarely occur in practice, since most problems come with additional constraints such as vehicle capacities and complex operating costs that are not only limited to distance.

 The VRP was first formally introduced in \cite{dantzig1959truck} and involves minimizing the total distance traveled by a fleet of vehicles serving a set of customers.  For a basic overview of the numerous VRP variants, we refer the reader to \cite{toth2002vehicle}. 
Solutions to industrial-sized VRPs are generally limited to heuristics built upon TSPs such as \cite{helsgaun2017extension} or more robust methods such as taboo searches like \cite{2001Taillard} and large neighborhood searches like \cite{pickup_and_delivery}. 
 Our work is concerned with a common variant of the VRP called the \struc{open vehicle routing problem (OVRP)} where vehicles are not forced to return to the depot after serving all the customers \cite{vincent2016open}. 
 The OVRP is very common in inbound logistics of the automotive industry, for example. 
 We define a specific type of OVRP that we consider throughout this article.
 It is an \struc{integer program (IP)}. 
\subsubsection{Open Vehicle Routing Problem with Stop Constraints and a Fixed Fleet}
\label{sec:OVRP}
The \struc{OVRP with stop constraints and a fixed number of vehicles} is a variant of the IP introduced in \cite{vincent2016open} like that in \cite{borcinova2017two} with the well-known extension of the \struc{subtour elimination constraints} of Miller-Tucker-Zemlin (MTZ) for the TSP \cite{miller1960integer}. 

In \cref{var_table,index_table}, we give an overview of the variables used in the OVRP formulation.

\begin{table}[h]
\begin{centering}
\begin{tabular}{ |c|c|c| } 
\hline
\textbf{Index} & \textbf{Sets} & \textbf{Description}\\ 
\hline
$i,j$ & $V$ & Generic set of all nodes including depot \\  
$0$ & -- & Depot node \\ 
$k$ & $K$ & Set of available vehicles \\ 
$n$ & $V\setminus\{0\}$ & Number of customer nodes \\ 
\hline
\end{tabular}
\caption{Definitions for the indices and labels appearing in the IP in \eqref{ovrpopt}-\eqref{xdomain}. }
\label{index_table}
\end{centering}
\end{table}
\begin{table}[h]
\begin{centering}
\begin{tabular}{| c | c | c |} 
\hline 
\textbf{Name} & \textbf{Domains} & \textbf{Description}\\ 
\hline 
$d_{ij}$ & $\mathbb{R}_{\geq 0}$ & Distance between nodes $i$ and $j$ \\ 
$C_k^{\textnormal{fixed}}$ & $\mathbb{R}_{\geq 0}$ & Fixed cost of each vehicle $k$\\ 
$C_k^{\textnormal{per dist}}$ & $\mathbb{R}_{\geq 0}$ & Cost per distance for vehicle $k$ \\ 
$C_k^{\textnormal{per stop}}$ & $\mathbb{R}_{\geq 0}$ & Cost for a vehicle $k$ to stop at a node in $V\setminus{\{0\}}$ \\ 
$\text{maxstop}_k$ & $\mathbb{Z}_{\geq 0}$ & Max number of stops allowed in $V\setminus{\{0\}}$ for each $k$  \\ 
$x_{ijk}$ & $\{0,1\}$ & Decision variable for if $k$ connects $i$ to $j$ \\ 
$u_{i}$ &  $\mathbb{Z}_{\geq 0}$ &  The order in which node $i$ is visited \\ 
$Z_{\text{OVRP}}$ &  $\mathbb{R}_{\geq 0}$ &  The objective value, viz. the total cost incurred \\ 
\hline 
\end{tabular}
\caption{Definitions for the cost factors, constraint bounds, and variables appearing in the IP in \eqref{ovrpopt}-\eqref{xdomain}.
}
\label{var_table}
\end{centering}
\end{table}
Let
\begin{equation}\label{totalcost}
    Z_{OVRP} = \sum_{k\in K}\sum_{i,j \in V}C_k^{\textnormal{per dist}} d_{ij}x_{ijk} +\sum_{k\in K}\sum_{\substack{i,j \in V\\ j \neq 0}}C_k^{\textnormal{per stop}} x_{ijk} + \sum_{k\in K}\sum_{j \in V\setminus\{0\}}C_k^{\textnormal{fixed}} x_{0jk}.
\end{equation}
Then the IP for the OVRP with stop constraints can be formulated as 
\begin{alignat}{2}
\label{ovrpopt}\min \qquad &Z_{OVRP}&&\\
\notag\\
\textnormal{subject to}\qquad\label{each_cust_visited_once} &\sum_{\substack{i\in V  \\k\in K\\ i\neq j}}x_{ijk} = 1, &&\text{ for all } j \in V\setminus \{0\},\\
\label{each_cust_left_once} &\sum_{\substack{j\in V\setminus \{0\}  \\k\in K\\ i\neq j}}x_{ijk} \leq 1,&&\text{ for all } i \in V\setminus \{0\},\\
\label{flow_constraint} &\sum_{\substack{i \in V\\ i\neq j}}x_{ijk}  -  \sum_{\substack{i\in V\setminus\{0\}\\i \neq j}}x_{jik} \geq 0, \qquad &&\text{ for all } k\in K, \; \text{ for all } j\in V\setminus\{0\},\\
\label{capcut1} &u_{i}-u_{j} + 1 \leq n(1-\sum_{\substack{k \in K}}x_{ijk}), \qquad &&\text{ for all } i,j \in V\setminus\{0\} \;\textnormal{with}\; i\neq j,\\
\label{maxstop} &\sum_{\substack{i\in V \\j\in V\setminus \{0\} \\i\neq j}}x_{ijk} \leq \text{maxstop}_k,\qquad &&\text{ for all } k \in K,\\
&\label{capcut2} 1 \leq u_{j} \leq n, \qquad  &&\text{ for all } j \in V\setminus\{0\},\\
\label{xdomain} &x_{ijk} \in \{0,1\}, \qquad  &&\text{ for all } i,j \in V, \; \text{ for all } k \in K.
\end{alignat}
Constraint \eqref{each_cust_visited_once} guarantees that every customer node is entered exactly once. Constraint \eqref{each_cust_left_once} guarantees that each customer is left \textit{at most} once.  Constraint \eqref{flow_constraint} is a type of flow constraint making sure that customer nodes cannot act as sources. Constraint \eqref{capcut1} is the MTZ-Subtour elimination constraint, forcing vehicles to have a tour connected with the depot. Constraint \eqref{maxstop} restricts the number of stops each vehicle can make. Constraints \eqref{capcut2} and \eqref{xdomain} define the domains for each of the variables.

 \subsubsection{The Instances}
 \label{sec:OVRP_instances}
In what follows, we describe two small real-world instances for OVRP.
We point out that our OVRP instances are small by industry standards. As a comparison, most OVRPs range from 50-200 customers.

\textbf{Data Preparation.}
The original real-world data includes geographic coordinates and other sensitive information that we anonymize: 
We replace the coordinates with a distance matrix using the great circle distance as well as remove names of locations, references to company names, vehicle makes and models. 
 
\textbf{Assumptions and Key Features.}
We briefly enumerate the underlying assumptions and key features of the two instances:
\begin{itemize}
    \item We label the instances by \texttt{OVRP15} and \texttt{OVRP11}, with 14 and 10 customer nodes, respectively. Each instance has a single depot node. Note that this leads to 15 and 11 total nodes, respectively.
    \item The instances are OVRPs with a homogeneous fleet of vehicles. This means that $C_k^{\textnormal{fixed}}$, $C_k^{\textnormal{per dist}}$, and $C_k^{\textnormal{per stop}}$ in \cref{var_table} are all independent of $k$. 
    \item There are no capacity constraints, since the total demand in the original customer problem did not exceed a single vehicle capacity.
    \item Each vehicle cannot make more than three stops, i.e., $\textnormal{\text{maxstop}}_k = 3$ for all $k$ in \eqref{maxstop}.
    \item We assume that each instance uses the minimum number of vehicles needed to solve the problem, which comes from the maxstop constraint in Equation \eqref{maxstop}: 
    \begin{align*}
        &\textnormal{mininum number of vehicles for } \texttt{OVRP15} = \left \lceil\frac{14}{\textnormal{3}}\right \rceil= 5,\\
        &\textnormal{mininum number of vehicles for } \texttt{OVRP11} = \left\lceil \frac{10}{\textnormal{3}}\right\rceil = 4.
    \end{align*}
    We constrain the number of vehicles to the smallest possible number in order to keep the QUBOs as small as possible. 
    Notice that under this assumption, the fixed costs in Equation \eqref{totalcost} no longer play a role.
\end{itemize}

\subsection{Slotting Problems}
\label{sec:slotting_problems}
The field of logistics offers a variety of non-linear optimization problems. Of particular interest for this work are problems that are inherently quadratic optimization problems. 
For example, optimally storing and retrieving items within warehouses is related to quadratic assignment problems \cite{ccela1999quadratic} and routing problems, such as the TSP. We introduce a type of slotting problem, which is an \struc{integer quadratic program (IQP)}.

\subsubsection{Affinity-based Slotting}
\label{sec:asp}
 We note that the IQP discussed in this section is essentially a special case of a \struc{generalized quadratic assignment problem (GQAP)} \cite{lee2004generalized}. Many works have studied GQAPs in the context of warehouses where instead of flow matrices, affinity matrices are used in the objective function \cite{Li2021, Chuang2012, Jiang2019, Lee2020, Li2016,HTX2023}. 
  Warehouses aim to ease the storing and procurement of goods for a company. Each warehouse stores an amount of products or materials of different \struc{material types $\mathcal{M}$} available for purchase or further shipment.
  Based on order history, one can determine which material types are generally ordered together and thus also need to be retrieved or ``picked" together. 
  This is the concept of \struc{pairwise affinity}. 
  We represent pairwise affinity for a pair of material types $m,n \in \mathcal{M}$ by an element of an \struc{affinity matrix} $\struc{\sigma}\in \mathbb{R}_{\geq 0}^{|\mathcal{M}|\times |\mathcal{M}|}$. 
  An element of $\sigma$ is generally larger for materials that are often ordered together, and smaller otherwise. 
  There are numerous ways to define $\sigma$ in the literature such as the \struc{Jaccard Measure} \cite{jaccard1901distribution} and the \struc{Bindi Measure} \cite{bindi2009similarity}, see \cite{phdthesis} for an overview.
  
  Given $\sigma$ and a set of material types $\mathcal{M}$, assume that we have access to a sequence of orders, and a set of storage locations that are partitioned into aisles $A_1,...,A_k$.  
  Each aisle, for $j=1,\dots,k$, has a \struc{capacity $|A_j|$}, which corresponds to the amount of material types that can fit in a single aisle. 
  The goal is to assign each material type $m\in\mathcal{M}$ to an aisle $A_j$
  such that we minimize the \struc{cost} of retrieving all orders.

  In this work, we consider the \struc{affinity-based slotting problem (ASP)},
  where the cost function aims to minimize the number of \struc{aisle changes} in order to find the minimal travel distance for retrieving all materials. 
   This problem 
   can be shown to be NP-hard, for example, by reducing the problem to a perfectly balanced bipartition problem \cite{JulianThesis,garey1974some}. 
   To model the ASP, let $x_{mj}\in \{0,1\}$ be equal to 1 if material type $m$ is assigned to aisle $A_j$, and 0 otherwise. 
   If different material types $m \neq n$ are placed in different aisles $A_i$ and $A_j$ with $A_i\neq A_j$, we impose a price proportional to their affinity $\sigma_{mn}$. 
   We formulate the ASP as an IQP in \eqref{ASP-1}.
\begin{alignat}{3}
&\min Z_{ASP}= &&\min \sum_{m,n\in\mathcal{M}}\sum_{\underset{i\not = j}{i,j=1}}^k\sigma_{mn}x_{mi}x_{nj}&&\label{ASP-1}\\
&\textnormal{subject to}&&\sum_{j=1}^kx_{mj}=1 && \text{for all }m\in\mathcal{M}\label{exactlyoneaisle}\\
&&&\sum_{m\in\mathcal{M}}x_{mj}\leq |A_j|&&\text{for all }j=1,...,k\label{aislecapcity}\\
&&&x_{mj}\in\{0,1\} &&\text{for all }m\in \mathcal{M},\; j=1,\dots,k. \label{aspvariables}
\end{alignat}
 
The constraint \eqref{exactlyoneaisle} forces each material to be assigned to exactly one aisle. Constraint \eqref{aislecapcity} ensures that the capacity of each aisle is not exceeded. The final constraint in \eqref{aspvariables} defines the domain for the problem variables. 
For more details, see \cref{appendix:details_problem_instances_ASP}. 

\subsubsection{The Instances}
\label{sec:asp_instances}
Given the ASP formulation in \eqref{ASP-1}, we provide three data sets using an affinity matrix $\sigma$ derived from real-world data.

\textbf{Data Preparation.}
These instances are \textit{based} on real order data. The order data comes from a warehouse containing over 4000 orders over a year-long period. 
The 30, 60, 90 most popular material types across all orders were determined, and then the affinity matrices $\sigma$ were constructed using the Jaccard measure \cite{jaccard1901distribution}.

\textbf{Assumptions and Key Features.}
We briefly describe the underlying assumptions and key features of the instances.  
We provide three instances that we label by \texttt{ASP30},  \texttt{ASP60},  \texttt{ASP90}, where the corresponding number denotes the number of material types. 
Each dataset has three aisles $A_j$ for $j=1,2,3$ with equal capacities $|A_1| = |A_2| = |A_3|$. The sum of the three aisle capacities equals the total number of materials. 
For example, \texttt{ASP30} is a dataset with a total capacity of 30 material types and aisle capacities of ten storage locations. 
Note that, we limit the number of aisles to  $k = 3$ in order to keep the problem sizes small enough to provide an optimal solution via a standard quadratic solver 
that we can compare with the SDP relaxations via SOS-SDP and HU.
However, limiting $k=3$ is an extreme underestimate of the real-world situation. Customers regularly have warehouses with $k>10$.

\subsection{Reformulation as QUBO}
\label{sec:QUBO_reformulation}
After formally defining the problems as IP and IQP models, we transform the formulation into a quadratic unconstrained binary optimization (QUBO) problem, see Equation \eqref{eq:QUBO}, which we further relax in the upcoming sections.
For the conversion from IP to QUBO, we follow the methodology developed by Glover et al. \cite{fred2019}. The same approach can be applied to convert an IQP to QUBO, as the constraints are linear and the objective function is quadratic.

Essentially, converting an IP/IQP model into a QUBO problem involves three steps:
\begin{enumerate}
    \item \textbf{Unconstraining}, which involves transforming constraints into terms that can be incorporated into the objective function. 
    \item \textbf{Binarization} is applied where integer variables are redefined as binary, ensuring the optimization problem fits within the binary framework of QUBO.
    \item \textbf{Penalization} to ensure that any violation of the original constraints is appropriately penalized.
\end{enumerate}
In \cref{sec:unconstraining,sec:binarization,sec:penalization} we show the conversion of the IQP formulation to its QUBO formulation in detail, using \texttt{ASP30} as an example.

The modeling and conversion process is carried out using \struc{\textsc{GAMS}} (\struc{General Algebraic Modeling System}), a platform for formulating and solving large-scale mathematical optimization problems. 
\textsc{GAMS} provides the flexibility to implement the transition from IP/IQP to QUBO efficiently and the computational power necessary to solve both types of problems. 
In particular, the \textsc{GAMS} reformulation tool \textsc{QUBO\_SOLVE}, see \cref{sec:ex_setup},
is employed to transform the original IP/IQP models in the \textsc{GAMS} language into the QUBO formulation \eqref{eq:QUBO}, facilitating the solving of the QUBO problem as described in the later sections.
This tool also supports solving the problem using either classical solvers such as e.g.~\struc{\textsc{CPLEX}} \cite{cplex} and \struc{\textsc{GUROBI}} \cite{gurobi} or hybrid quantum methods such as D-Wave's quantum annealing procedure \cite{DWaveLeap}. 
It is important to note that the objective value obtained after solving the QUBO to optimality will be equal to the optimal value of the original IP/IQP formulation, as both the models are equivalent, see \cite{AdiabaticQuantumComputingandQuantumAnnealing} and \cite{chang2022hybridquantumclassicalcomputing}.

\section{SOS-SDP Relaxation}
\label{sec:SOSSDP}
In this section we interpret QUBOs as a \struc{polynomial optimization problem on the boolean hypercube}.
We recall that in the QUBO formulation \eqref{eq:QUBO} we are given a symmetric matrix $Q~=~(q_{ij})_{1\leq i,j\leq n} \in\mathbb{R}^{n\times n}$, and our goal is to minimize $\Vector{x}^T Q \Vector{x}$ with $\Vector{x}\in\{0,1\}^n$. 
Rewriting this problem as a polynomial optimization problem, we aim to apply a \struc{sums of squares (SOS) hierarchy relaxation} that leads to solving a semidefinite program (SDP).
In \cref{section:polyopt} we give an overview on the interpretation as polynomial optimization problem and in \cref{section:sos_sdp_relaxations} we motivate how to derive the SDP relaxation from the SOS approach. 
We explain the software that we use for the SOS relaxation in \cref{section:tubs:software}.
\subsection{The QUBO as Polynomial Optimization Problem}
\label{section:polyopt}
Let $\mathbb{R}[\Vector{x}]= \mathbb{R}[x_1\; \hdots\; x_n]$ be the polynomial ring of \struc{real $n$-variate polynomials}.
Rewriting the target function from the given QUBO as in Equation \eqref{eq:QUBO} we can consider the following polynomial $p(\Vector{x}) \in \mathbb{R}[\Vector{x}]$ with $n$ variables of degree two: 
\begin{equation}
\begin{aligned}
    p(\Vector{x}) &= \Vector{x}^TQ\Vector{x}\label{TUBS:qubopoly}\\ &= \begin{bmatrix}
        x_1&x_2&\hdots&x_n
    \end{bmatrix} 
    \begin{bmatrix}
        q_{11} & q_{12} & \hdots & q_{1n}\\
        q_{12} & q_{22} & \hdots & q_{2n} \\
        \vdots & \vdots &\ddots & \vdots \\
        q_{1n} & q_{2n} & \hdots & q_{nn}
    \end{bmatrix}
    \begin{bmatrix}
        x_1\\x_2\\ \vdots \\ x_n
    \end{bmatrix} \\
    &= q_{11}x_1^2 + 2q_{12}x_1x_2 + \hdots + 2q_{1n}x_1x_n + q_{22}x_2^2 + \hdots + 2q_{2n}x_2x_n + \hdots  + q_{nn}x_{n}^2.
\end{aligned}
\end{equation}
Furthermore, we rewrite the binary constraints $\Vector{x} \in \{0,1\}^n$ as polynomial equality constraints:
\begin{equation}
    x_i(x_i - 1) = 0 \text{ for all } 1\leq i \leq n.\label{TUBS:eqconstraints}
\end{equation}
Using the polynomial target function from \eqref{TUBS:qubopoly} and the polynomial constraints \eqref{TUBS:eqconstraints} we obtain a \struc{constrained polynomial optimization problem (CPOP)} as interpretation of the given QUBO.
Equivalently to minimizing $p(\Vector{x})$, we can formulate this as a problem of \struc{nonnegativity}:
\begin{equation}
\begin{aligned}
			&\text{max}							
			& &	\lambda \in \mathbb{R}\\
			&\text{subject to}	&	&	p(\Vector{x}) - \lambda \geq 0 \text{ for all } \Vector{x}\in \{0,1\}^n.\\
		\end{aligned}\label{TUBS:nonnegativity}
\end{equation} 
Polynomial optimization problems are NP-hard in general.
However, the idea of the nonnegativity formulation in \eqref{TUBS:nonnegativity} offers the use of relaxations in terms of \struc{certificates of nonnegativity}.
These are sufficient conditions showing a polynomial $f(\Vector{x})$ is nonnegative, i.e., in the unconstrained case $f(\Vector{x})\geq 0$ for $\Vector{x}\in \mathbb{R}^n$, and the conditions are effectively computable.
A widely used certificate of nonnegativity is to show that a polynomial $f(\Vector{x})$ can be written as $f(\Vector{x}) = \sum_{i=1}^{\ell} s_i(\Vector{x})^2 \in \mathbb{R}[\Vector{x}]$, i.e., $f(\Vector{x})$ is a sum of squared polynomials $s_i(\Vector{x})$.
We refer to this relaxation as \struc{sum of squares (SOS)} relaxation and to a polynomial that can be written as such decomposition likewise as \struc{SOS} (polynomial).
For an overview on polynomial optimization and SOS see e.g.~\cite{Laurent:Survey, Lasserre:IntroductionPolynomialandSemiAlgebraicOptimization, Theobald:Book:RealAlgGeom}
and \cref{sec:nonnegativity_certificates} for an overview on other certificates of nonnegativity beyond SOS. 
\subsection{Deriving the SDP relaxation}
\label{section:sos_sdp_relaxations}

We aim to apply an SOS relaxation on \eqref{TUBS:nonnegativity} which leads to solving an SDP.
For an intuition, we first consider unconstrained polynomial optimization. 

\textbf{Gram Matrix Method.}
Let $f(\Vector{x})$ be a polynomial in $n$ variables and of degree $2d$. We want to determine whether $f$ is SOS.
Let $\Vector{m_d}$ be the vector that contains all monomials (e.g.~$1, x_1, x_2, x_1^2, x_1x_2,...$) with $n$ variables up to degree $d$.
A symmetric matrix $G$ is positive semidefinite if and only if there exists a \struc{Cholesky decomposition}, i.e., there exists a matrix $B$ with $B^TB = G$. 

We aim to find a \struc{Gram matrix} $G \in \mathbb{R}^{{\binom{n+d}{d}}\times {\binom{n+d}{d}}}$ containing the corresponding coefficients of $f(\Vector{x})$ such that $\Vector{m_d}^T\;G\;\Vector{m_d} = f(\Vector{x})$. If and only if $G$ is positive semidefinite and thus permits the Cholesky decomposition $G=B^TB$, then 
\begin{align*}
    f(\Vector{x}) = \Vector{m_d}^T\;G\;\Vector{m_d} = \Vector{m_d}^TB^TB\;\Vector{m_d} = (B\;\Vector{m_d})^TB\;\Vector{m_d} = \sum s_i(\Vector{x})^2,
\end{align*}
i.e., we can write $f(\Vector{x})$ as a sum of squared polynomials $s_i(\Vector{x})$. 

\textbf{Unconstrained Polynomial Optimization.}
The Gram-Matrix-Method can be combined with semidefinite programming  \cite{parrilo2000structured}:
Using the definition of $\Vector{m_d}$ as above, the unconstrained problem of nonnegativity $\max_{\lambda\in \mathbb{R}}\{p(\Vector{x})-\lambda\ \geq 0 \text{ for all } \Vector{x}\in\mathbb{R}^n\}$ turns into the SDP relaxation:
\begin{equation*}
\begin{aligned}
			&\text{max}							
			& &	\lambda \in \mathbb{R}\\
			&\text{subject to}	&	&	p(\Vector{x}) - \lambda = \Vector{m_d}^T\;G\;\Vector{m_d} \\
            & & & G \succcurlyeq 0.
		\end{aligned}\label{TUBS:unconstrained_SOS}
\end{equation*} 
We optimize over the entries of a Gram matrix $G$ with linear constraints that refer to the coefficients of the polynomial $p(\Vector{x})-\lambda$ with the requirement that $G$ is positive semidefinite. 
Thus, we have to solve an SDP.

\textbf{Constrained Polynomial Optimization.} 
After this intuition in the unconstrained case we again consider our original CPOP in \eqref{TUBS:nonnegativity} which we derived from the given QUBO.
Putinar's \cite{putinar} Positivstellensatz uses the idea to express a polynomial in terms of polynomial constraints $h_i(\Vector{x}) = 0,$ for $ i=1,\dots, n$ and SOS polynomials $s_i(\Vector{x})$ for $ i=0,\dots,n$ to certify nonnegativity. This leads to a method known as \struc{Lasserre's hierarchy} \cite{lasserre} that tackles the problem:
\begin{equation}
\begin{aligned}
			&\text{max}							
			& &	\lambda \in \mathbb{R}&\\
			&\text{subject to}	&	&	p(\Vector{x}) - \lambda = s_0(\Vector{x}) + \sum_{i=1}^n s_{i}(\Vector{x})h_i(\Vector{x})&\\
			&&&	\text{degree}(s_0(\Vector{x})) \leq 2d \text{ and }\text{degree}(s_i(\Vector{x})h_i(\Vector{x})) \leq 2d &\text{for all } 1\leq i\leq n\\
            &  &  & s_0(\Vector{x}) \text{ is SOS and } s_i(\Vector{x}) \text{ are SOS} &\quad\text{for all } 1\leq i\leq n
		\end{aligned}\label{TUBS:lasserre}
\end{equation}
Problem \eqref{TUBS:lasserre} again results in solving an SDP: We aim to find polynomials $s_i(\Vector{x})$ that are SOS.
The polynomials $s_i(\Vector{x})$ are unknown and we optimize over their corresponding Gram matrices, restricting their degree.
By increasing the allowed degree $2d$ we compute a sequence of SDP relaxations. 
This sequence leads to \struc{sequence of lower bounds} that -- under mild assumptions, are fulfilled in the QUBO case -- eventually converges to the actual optimum of \eqref{TUBS:nonnegativity} \cite{Fawzi2016, lasserre}. 
In fact, the approach by Lasserre \cite{lasserre} uses the theory of \struc{moments} and is dual to \eqref{TUBS:lasserre}. 
For further details, see \cite{Lasserre:IntroductionPolynomialandSemiAlgebraicOptimization,Theobald:Book:RealAlgGeom}.

Note that, the involved Gram matrices in the sequences of lower bounds grow very large in practice. They have a size of $\binom{n+2d}{2d}\times\binom{n+2d}{2d}$ for a polynomial with $n$ variables and degree $2d$. The size increases exponentially, as we increase the degree in the hierarchy.
For our problems, the number of variables $n$ ranges between $96\leq n\leq 3380$ and we see that for our problems, even the first order Lasserre relaxation with $2d = 2$ is challenging for SDP solvers, see \cref{sec:experiments}.

\subsection{Solving the SOS relaxation}
\label{section:tubs:software}
In order to compute the SOS-SDP relaxation as described before we work with academic software used in polynomial optimization to derive the SDP relaxations of the QUBOs.
Then an SDP solver is called internally to solve those SDPs.
We give give an overview on the tested software below.
In preliminary experiments, we additionally tested the software \textsc{gloptipoly} \cite{gloptipoly} on \textsc{matlab}. We did not pursue this further since the software struggled with large input matrices.

\textbf{\textsc{SumOfSquares.jl}}
The \textsc{julia} package \struc{\textsc{SumOfSquares.jl} (\textsc{SOS.jl})} \cite{SOS_julia2017, SOS_julia2019} performs unconstrained and constrained SOS optimization as described above.
We input the given QUBO matrix and convert it to a polynomial representation.
The \textsc{SOS.jl} package allows to input the relaxation order as in \eqref{TUBS:lasserre} automatically derives the corresponding SDP.
Then, internally, the package calls an SDP solver to solve the resulting SDP.
We test the internal SDP solvers \struc{\textsc{Mosek}} \cite{mosek} and \struc{\textsc{Loraine}} \cite{loraine2023}.

\textbf{\textsc{TSSOS}}  
The \textsc{julia} software \struc{\textsc{TSSOS}} \cite{TSSOS}\cite{TSSOSDok} by Wang, Magron and Lasserre improves the above described moment SOS hierarchy along with tools to exploit sparsity in the given polynomial \cite{wang2020tssos}\cite{wang2020chordaltssos} (``TS'' stands for ``term sparsity'') by exploring block structures in the matrices. 
Similar to \textsc{SOS.jl}, \textsc{TSSOS} uses the polynomial representation and then automatically derives the corresponding SDPs.
As SDP solver, internally, \textsc{TSSOS} provides an interface for \textsc{Mosek} and \struc{\textsc{COSMO}} \cite{cosmo}.
Especially, the extension \textsc{CS-TSSOS} \cite{wang2021cstssos} (``CS'' stands for ``chordal sparsity'') is suitable for large-scale polynomial optimization problems.
Earlier works on chordal sparsity include the software \textsc{SparsePOP} \cite{sparsePOP}, that we do not test here.

\textbf{SDP solvers} Inside of \textsc{SOS.jl} and \textsc{TSSOS}, we call the following SDP solvers as described above: The SDP solver \textsc{Mosek} is a state-of-the-art but general purpose interior-point method SDP solver. 
As shown in \cref{sec:experiments}, while we achieve the best results on our instances with \textsc{Mosek} inside \textsc{TSSOS}, we experience that the large-scale  SDPs cannot be solved. 
Additionally, we use the low-rank interior-point SDP solver \textsc{Loraine} for \textsc{julia}.
It has been reported that often with \textsc{Loraine} the second order relaxation already converges to the optimal solution, outperforming \textsc{Mosek} on randomly generated MaxCut and QUBO instances \cite{habibi:hal-04076510}.
The solver \textsc{COSMO} does not perform interior-point methods but is based on the alternate direction of multipliers method (ADMM) and is suitable for large-scale convex optimization problems.

Next to the polynomial optimization software \textsc{TSSOS} and \textsc{SOS.jl}, we also solve the SOS-SDP relaxation 
directly via available academic SDP solvers:
The \textsc{matlab} solver \struc{\textsc{ManiSDP}} \cite{wang2024solvinglowranksemidefiniteprograms} based on \textsc{Manopt} \cite{manopt} is tailored for low-rank SDPs and uses manifold optimization. 
It has been used to solve the second order Lasserre relaxation ($2d = 4$) on QUBOs up to $n=120$ variables, outperforming \textsc{Mosek}. \cite{wang2024solvinglowranksemidefiniteprograms}.
For this solver, we use the Ising formulation, see \cref{subsection:binary_conversion}, as input. 
Likewise, we use the \struc{\textsc{Splitting Conic Solver}}\struc{(\textsc{SCS})} \cite{scs} with its default settings, as a standard an open source convex cone solver, directly on this formulation.

\section{Hamiltonian Updates: Classical and Quantum}
\label{sec:HU}
In this section, we give an overview on the SDP relaxation via the \struc{Hamiltonian Updates algorithm (HU)}, developed in \cite{Brandao2022fasterquantum} and improved by \cite{henze2025}. 
The advantage of HU is its favorable scaling in the problem size compared to standard SDP solvers at the cost of a rather expensive polynomial dependence in the precision. 
This makes HU more suitable to applications where high precision is less important:
This is for example the case when we use the solution of the SDP as input to a randomized rounding procedure, which is approximate by its nature.
Additionally, through the randomized rounding procedure, HU provides lower and upper bounds for the optimal value, as well as a feasible solution for the original problem.
\subsection{Problem formulation}
For HU, we consider a slightly reformulated QUBO problem of the form 
\begin{equation}
\begin{aligned} \label{eq:ising}
   &\max_{\Vector{x}\in\mathbb{R}^n} & &\sum_{i,j=1}^{n} C_{ij}x_i x_j\\
   &\text{subject to } & & x_i=\pm1,
\end{aligned}
\end{equation}
which we call an \struc{Ising model} with a cost matrix $C\in\mathbb{R}^{n\times n}$.
As shown in \cref{subsection:binary_conversion}, one can easily convert a problem given in QUBO form \eqref{eq:QUBO} into the Ising formulation \eqref{eq:ising}. 
We can then relax \eqref{eq:ising} into the following SDP:
\begin{equation}
    \begin{aligned} 
        &\max_{X\in\mathbb{R}^{n\times n}} & & \langle C,X\rangle\\\label{eq:sdp_relaxation}
        &\text{subject to } & & \diag(X)=1\\
        & & &  X \succcurlyeq 0.       
    \end{aligned}
\end{equation}
Goemans and Williamson \cite{Goeman_Williamson} developed a \struc{randomized rounding} procedure that,
given a solution for the SDP relaxation \eqref{eq:sdp_relaxation}, 
outputs a feasible solution for the original Ising problem. 
In addition, for cost matrices of specific forms, this method comes with rigorous approximation guarantees, e.g.~applied to the original MaxCut formulation with nonnegative edge weights, the randomized rounding procedure by Goemans and Williamson gives the famous result of an expected value that is at least $0.87856$ times as large as the optimal value \cite{Goeman_Williamson}.

For the randomized rounding, we first
factorize the optimizer $X^*\in \mathbb{R}^{n\times n}$ of the SDP in (\ref{eq:sdp_relaxation})
as $X^*~=~B^T B$ with a symmetric matrix $B\in\mathbb{R}^{n\times n}$. 
Then the sign of the column-wise projection of $B$ onto a Gaussian random vector $\Vector{g}_i\in\mathbb{R}^n$ gives the entries of a vector $\Vector{x}\in\{\pm1\}^n$ that is an approximate solution to \eqref{eq:ising}. 
We summarize this rounding procedure as
\begin{equation}\label{eq:round1}
    \begin{aligned}
        & \text{compute} \quad && B=\sqrt{X^*}, \quad &  \\
    & \text{sample} \quad && {\Vector{g}_i \overset{\mathrm{idd.}}{\sim} \N(0,1),} &\quad i=1,\dots,n
    \\
    & \text{compute} \quad && {x_j \gets \sgn\left(\sum_{i}\Vector{g}_i B_{ij}\right),} &\quad j=1,\dots,n.
    \end{aligned}
\end{equation}
where $\sqrt{X^*}\in \mathbb{R}^{n\times n}$ denotes a symmetric matrix with  $\sqrt{X^*}^T\sqrt{X^*} = X^*$ and $\N(0,1)$ is a normal distribution with a mean $0$ and variance $1$.

\subsection{Solving the SDP using Hamiltonian Updates}

The idea behind HU is to express $X$ as a \struc{Gibbs state}, an object known from quantum mechanics. The physical context of this is not further relevant here since the HU routine itself is purely classical. Nevertheless, this inspiration from quantum mechanics is a natural choice for positive semidefinite matrices, and gives rise to a possible quantum speedup \cite{vanapeldoorn2018}.
The Gibbs state is given by
\begin{align} \label{eq:gibbs}
	\rho_H = \frac{\exp(-H)}{\tr(\exp(-H))} \in \mathbb{R}^{n\times n}
\end{align}
for some symmetric $H\in\mathbb{R}^{n\times n}$, called the \struc{Hamiltonian}, and we set $X=n\rho_H$. 
This construction makes $X$ manifestly positive semidefinite. 
Then, by introducing a precision parameter $\epsilon\geq 0$, one can reformulate the SDP in \eqref{eq:sdp_relaxation} as feasibility problem with \struc{threshold objective value} $\struc{\gamma}\in \mathbb{R}$:
\begin{equation}
\begin{aligned}
    &\mathrm{find} &&H\in \mathbb{R}^{n\times n} \\
    &\text{subject to } & &\tr(C \rho_H) - \frac{\gamma}{n} < \epsilon \norm{C}\\
    & & &\sum_{i=1}^n \Big|(\rho_H)_{ii} - \frac 1 n\Big| < \epsilon,
\end{aligned}\label{eq:HU_feasibility_problem}
\end{equation}
where $\|\cdot\|$ denotes the operator norm.
We can approximately determine the minimal threshold objective value $\gamma^*\in [-n\norm{C}, n\norm{C}]$ of \eqref{eq:sdp_relaxation} via a binary search. 
Here, the task of HU is to either find an $H$ that fulfills both constraints in \eqref{eq:HU_feasibility_problem}, or prove that no feasible $H$ exists for a given value of $\gamma$. 
We can transform a feasible solution for the SDP relaxation \eqref{eq:HU_feasibility_problem} to a solution of the Ising problem \eqref{eq:ising} by applying randomized rounding as shown in \eqref{eq:round1}. This means, we retrieve a feasible solution vector to the original Ising formulation and additionally an upper bound to the optimal value.
Conversely, if $\gamma$ is infeasible for the SDP relaxation \eqref{eq:HU_feasibility_problem}, it inherently serves as a lower bound for the Ising problem \eqref{eq:ising}, as given from the nature of relaxations. 
We can convert both, upper and lower bounds, for the Ising problem to bounds for the QUBO \eqref{eq:QUBO}, and thus also for the original IP \eqref{totalcost} and IQP \eqref{ASP-1}, by reversing the steps of \cref{subsection:binary_conversion}.

The HU algorithm is an iterative process, where in each step one estimates the violation of the constraint \eqref{eq:HU_feasibility_problem} and updates $H$ accordingly. 
We detect infeasible instances based on principles from quantum statistical mechanics involving an upper bound on the \struc{relative entropy distance} between the initial Gibbs state $\rho_0=\Id_{n\times n}/n$ and a possible feasible solution state $\rho^*~\in~\mathbb{R}^{n\times n}$. 
To check the constraints, the matrix exponential $\exp(-H)$ needs to be computed in each iteration of HU, which is by far the most expensive part of the algorithm.
Here, \struc{quantum convex optimization} comes into play \cite{Brandao2017}.
The central idea is to read (\ref{eq:gibbs}) as a \struc{physical} procedure, where $H$ is a Hamiltonian operator and $\rho_H$ a thermal quantum state.
The problem of preparing such a Gibbs state given a classical description of the Hamiltonian is a well-studied problem in quantum computing \cite{vanApeldoorn2020quantumsdpsolvers}.
This idea can be turned into a quantum algorithm for HU, with running time
$\tilde\O(n^{1.5}s^{0.5}\epsilon^{-5})$,
compared to
$\tilde\O(\min\{n^3\epsilon^{-2}, n^2s\epsilon^{-3}\})$ 
for a classical implementation, where $n$ refers to the size of the matrix $H\in \mathbb{R}^{n\times n}$, $s$ denotes the maximal number of non-zero entries in the matrix and $\epsilon$ is the precision as defined in \eqref{eq:HU_feasibility_problem}.
For the details, see \cite{henze2025}.

\section{Experimental Results}
\label{sec:experiments}
In the upcoming \cref{sec:experiments}, we compare the SOS-SDP relaxation, see \cref{sec:SOSSDP}, and Hamiltonian Updates, see \cref{sec:HU}, computationally. 
We apply the SDP relaxations of both approaches to the QUBOs that we derive from the IP/IQP formulations of the problem instances from \cref{sec:OVRP,sec:asp}. 
As benchmark values for comparison, we provide the optimal values of the respective IPs/IQPs as well as the optimal/approximate values to the QUBO formulations by standard solvers, see \cref{sec:mip_solutions,sec:qubo_solver_results}.
\subsection{Experimental Setup}
\label{sec:ex_setup}
We give a technical overview on our hardware and software setup:

\begin{itemize}
   \item \textbf{Hardware and system information:} We use Ubuntu 18.04.6 LTS with an Intel(R) Core(TM) i7-8700 CPU @ 3.20GHz processor, six cores and two threads per core and 15GB of RAM by default for all computations. We compute selected examples using HU on an Nvidia Geforce RTX 4090 GPU.
   \item \textbf{Code and data availability:}
We provide our data and executable scripts via \url{https://moto.math.nat.tu-bs.de/appliedalgebra_public/sdp_qubo_benchmarking}.
\item \textbf{Software:}
We derive the QUBO formulations via the \textsc{GAMS} \textsc{QUBO\_SOLVE} tool available at \url{https://github.com/GAMS-dev/qubo}.
Additionally, we use \textsc{GAMS} to provide IP/IQP solutions and QUBO solutions via the standard solvers \textsc{CPLEX} \cite{cplex} and \textsc{GUROBI} \cite{gurobi}. 
For the SDP relaxations based on SOS we rely on several academic and commercial software \cite{TSSOS, SOS_julia2019, mosek, cosmo, loraine2023, wang2024solvinglowranksemidefiniteprograms, scs}. For Hamiltonian Updates, we use the implementation  introduced in \cite{henze2025} and available at \url{https://zenodo.org/records/14871936}.
\end{itemize}

\subsubsection{Instances}
\label{sec:experiments_instances}
As laid out in detail\footnote{In particular, the QUBOs we consider have a constant offset $c \in \mathbb{Z}$, depending on the penalty from the QUBO reformulation. 
Thus, for a symmetric matrix $Q\in \mathbb{R}^{n \times n}$, we tackle the problem
  $  \min_{\Vector{x} \in\{0,1\}^n} \Vector{x}^T Q \Vector{x} + c .$
This problem differs from our initial definition by the offset $c$.
However, since $c$ does not alter the minimizer, we first solve a QUBO of the form \eqref{eq:QUBO} and later add $c$ to the solution.} in \cref{sec:OVRP_instances,sec:asp_instances}, we consider two OVRP instances emerging from an IP that we denote by \texttt{OVRP11} and \texttt{OVRP15} respectively, as well as three ASP instances, emerging from an IQP, that we denote by  \texttt{ASP30},  \texttt{ASP60} and  \texttt{ASP90}.
For the latter instances, we use multiple different penalties in the QUBO formulation to construct multiple QUBOs in order to explore if and how the penalties influence the solution of the SDP.
We distinguish the different QUBOs emerging from the same IQP by a suffix to their name with the respective penalty.
We give an overview over the five instances in \cref{tab:overview_intances}.

	\begin{table}[h!]
		\begin{center}
			\begin{tabular}{|l||c|c|c||c|c|}
            \hline
				&	\multicolumn{3}{c||}{\textbf{ASP}} & \multicolumn{2}{c|}{\textbf{OVRP}} \\
                \hline
				Instance &	 \texttt{ASP30} &  \texttt{ASP60} &  \texttt{ASP90} &    \texttt{OVRP11} &    \texttt{OVRP15} \\
                \hline
				QUBO Matrix size &	$n=96$ & $n=186$ & $n=276$ & $n=1750$ & $n=3380$ \\
                \hline 
                Sparsity & $49.9\%$& $53.5\%$ &$56.3\%$ &$92.2\%$ &$93.5\%$\\
                \hline
				Number of QUBOs &	$14$ & $4$ & $3$ & $1$ & $1$ \\
				 \hline
			\end{tabular}
		\end{center}
        \caption{Overview on the five \textbf{problem instances}. We display the size of the matrix $Q\in \mathbb{R}^{n\times n}$ and the sparsity as percentage of zero elements in $Q$ after the QUBO reformulation. For the ASP instances, we each provide multiple QUBO formulations that differ by penalties used in the construction.}
        \label{tab:overview_intances}
	\end{table}
The OVRP instances are very large with matrix sizes $n=1750$ and $n=3380$ respectively, yet relatively sparse with each over $92\%$ of zero entries. However, this yields still 238350 non-zero elements for   \texttt{OVRP11} and 740864 non-zero elements for \texttt{OVRP15}.
Based on analyzing the singular values, the QUBO matrices of \texttt{OVRP11} and \texttt{OVRP15} are of relatively high rank, with numerical ranks of 1480 and 3002, respectively. 

The ASP instances are smaller with $n=\{96, 186, 276\}$ and all have full numerical rank.
Moreover, in each matrix about half of the entries are not zero.
Note that the penalties in the QUBO formulation do not change the number nor position of the non-zero elements.
Thus, the value for sparsity for each version of the QUBOs emerging from the ASP instances is the same.

\begin{figure}[h]
    \begin{subfigure}[b]{0.45\textwidth}
        \centering
        \includegraphics[width=1\linewidth]{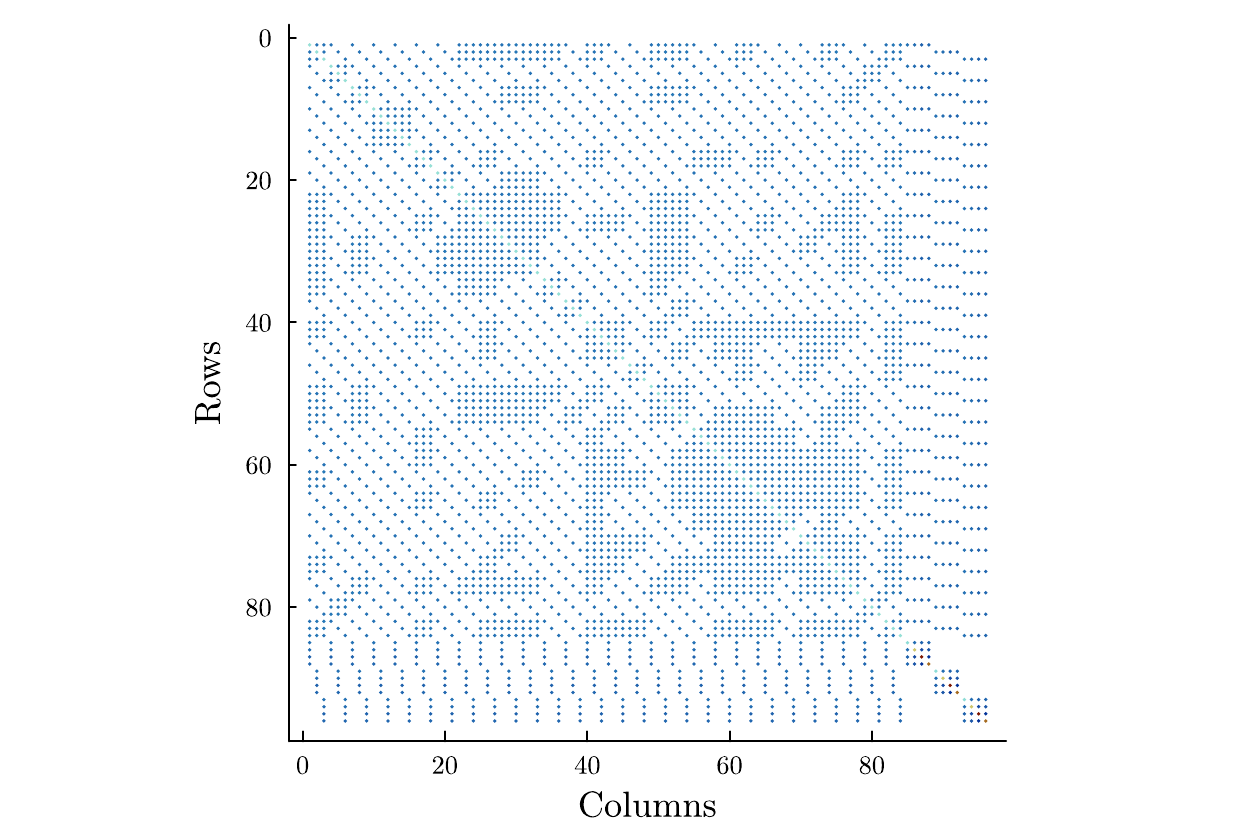} 
        \caption{Sparsity Pattern for \texttt{ASP30}-250, all dots resemble a non-zero entry.}
        \label{fig:image1}
    \end{subfigure}
    \hspace{0.8cm}
    \begin{subfigure}[b]{0.45\textwidth}
        \centering
        \includegraphics[width=\linewidth]{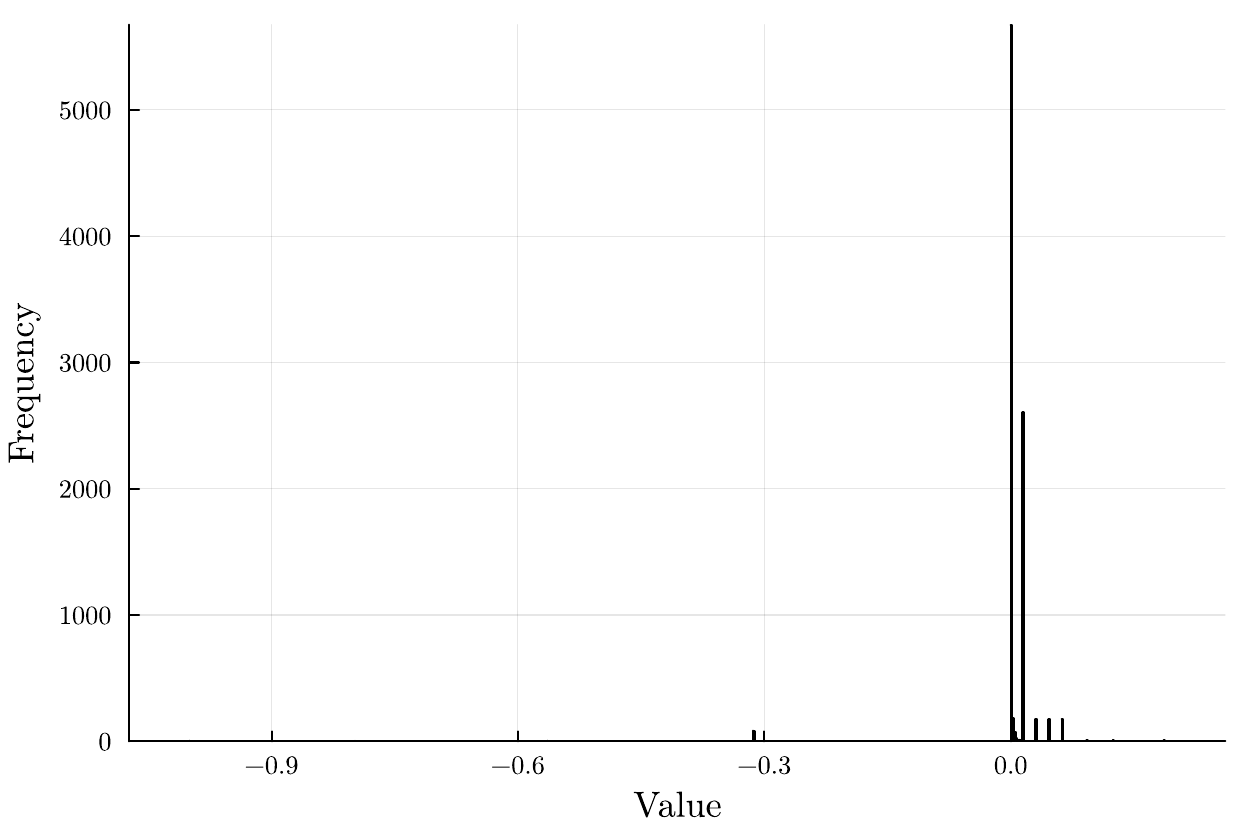} 
        \caption{Histogram of the distribution of matrix entries \texttt{ASP30}-250.}
        \label{fig:image2}
    \end{subfigure}
    \caption{\textbf{Sparsity} and \textbf{distribution of non-zero entries} of the matrix corresponding to QUBO \texttt{ASP30}-250. All values are scaled by dividing with the largest entry of absolute value. The largest values by absolute value lie on the diagonal, this is indicated by different colors on the diagonal.}
    \label{fig:plots_asp30p250}
\end{figure}

In \cref{fig:plots_asp30p250} we exemplarily show the sparsity pattern of the matrix entries for  \texttt{ASP30} with penalty 250, along with a histogram of the entries. It is visible that the entries are not normally distributed. 
In the display, all matrix entries are divided by the maximal absolute value.
\subsubsection{Attributes for comparison}
\label{sec:experiments_attributes}
In this section we introduce the different attributes that we use to compare the SOS-SDP methods and Hamiltonian Updates. 
While both approaches are SDP relaxations of the given QUBO they contain specific features that are not directly comparable.
With the goal in mind of benchmarking SDP relaxations for QUBO formulations from a practitioner's perspective on the specific given problem instances, we choose to compare in general two attributes that we attain from both approaches.

We label the \struc{optimal value} of the ASP in the IQP formulation in Equation \eqref{ASP-1} as \struc{$Z_{\text{ASP}}^*$} and the \struc{optimal value} of the OVRP in the IP formulation in Equation \eqref{totalcost} as \struc{$\zVRP$}. For simplicity, we refer to both optimal values of ASP and OVRP as \struc{$\Zopt$}.

\begin{itemize}
    \item \textbf{Lower bound:} Both SDP relaxations via SOS and HU provide a \struc{lower bound} \struc{\zSOS} and respectively a \struc{lower bound} \struc{\zHU} to the optimal value\footnote{Both relaxations provide a lower bound to the optimal QUBO value. 
    The QUBO is a reformulation of the corresponding IP/IQP with the same optimal value $\Zopt$. 
    The specific optimal values $\Zopt$ for our problem instances are positive and displayed in \cref{tab:MIP_results}.} $\Zopt$.
    That means, we have $\zSOS, \zHU \leq \Zopt$.
    In order to evaluate the quality of the lower bounds $\zSOS$ and $\zHU$ and compare them among each other, we state their distance to the optimal value $\Zopt$.
    We do this via the
    \struc{absolute differences} \struc{$\DeltaSOSabs$} and \struc{$\DeltaHUabs$}, denoted by
    \begin{align*}
        \DeltaSOSabs = |\zSOS-\Zopt| \text{ and } \DeltaHUabs = |\zHU-\Zopt|,
    \end{align*}
    and via the 
    \struc{relative differences \DeltaSOS} and \struc{\DeltaHU} to the optimal value $\Zopt$, denoted by
    \begin{align*}
        \DeltaSOS = \frac{\DeltaSOSabs}{\Zopt} \text{ and } \DeltaHU = \frac{\DeltaHUabs}{\Zopt}.
    \end{align*}
    
    \item \textbf{Running Time:} We record the running times \struc{\tSOS} and \struc{\tHU} respectively of each SDP relaxation that is spent to provide the above lower bounds $\zSOS$, $\zHU$.
\end{itemize}

Each approach provides multiple different parameters that influence the quality of the lower bounds as well as the running time. 
These parameters are individual to each approach. 
Likewise, we collect additional data that is not directly comparable to the other approach. 
We briefly describe these features for both methods:

\smallskip
\textbf{Additional parameters and tracked data for SOS-SDP:} 
\begin{itemize}
\item The SOS hierarchy as in \eqref{TUBS:lasserre} 
provides a sequence of lower bounds to the optimal value, depending on the degree $2d$ of the corresponding polynomials.
In our experiments, we perform first order relaxations ($2d=2$) and second order relaxations ($2d = 4$) if possible due to the size of the corresponding SDP. 
\item As laid out in \cref{section:tubs:software}, we use different polynomial optimization software as well as internal SDP solvers.

\item Note that a solution recovery for the SOS methods, i.e., attaining the vector $\Vector{x}^*\in \{0,1\}^n$, is only guaranteed if the optimal value is reached within the hierarchy, see e.g.~\cite{henrionlasserre, Laurent:Survey}. In our experiments we do not reach the optimal value. 
\end{itemize}

\smallskip
\textbf{Additional parameters and tracked data for HU:}
\begin{itemize}
    \item As described in \cref{sec:HU}, a crucial parameter in Hamiltonian Updates is the precision $\epsilon$ of the feasibility problem \eqref{eq:HU_feasibility_problem}. 
For a high precision, the running time increases accordingly.
In our experiments, we test values of $\epsilon$ in $\{10^{-2}, 10^{-3}, 10^{-4}\}$.
\item The quantum version of HU gives an approximate number of quantum CNOT gates that would be used to perform the algorithm on a quantum computer, see \cite{henze2025} for further details on the gate count. We assume the current best time of $6.5\times 10^{-9} s$ that could physically be realized to perform an isolated quantum gate operation \cite{Chew2022} and multiply by the approximate number of CNOT gates.
This gives an optimistic approximation of the running time \struc{\tHUQuantum} of the SDP relaxation via HU on a quantum computer.
\item The HU implementation provides the opportunity to be executed on a GPU via cupy \cite{cupy} if available instead of the default CPU. Execution on GPU is faster. 
\item 
The randomized rounding procedure of Hamiltonian Updates, see \eqref{eq:round1}, yields an approximate solution vector $\Vector{x}_{\text{HU}}\in \{0,1\}^n$ that provides an \struc{upper bound $\zHUupper$} to the optimal value $\Zopt$.
We define the absolute difference $\struc{\DeltaHUabsupper}$ and relative difference $\struc{\DeltaHUupper}$ to the optimal value $\Zopt$, respectively, as
\begin{align*}
    \DeltaHUabsupper = |\zHUupper - \Zopt| \text{ and } \DeltaHUupper = \frac{|\zHUupper - \Zopt|}{\Zopt}.
\end{align*}
\end{itemize}
We execute all SDP relaxations with the parameters as described above. 
Since we try to take a view point of a practitioner's perspective, we generally enforce time limits for the computations and stop them accordingly, see \cref{tab:timelimits} in the appendix for details.

\subsection{Benchmarks from IP/IQP and QUBO formulations}
\label{sec:optimal_mip_qubo_results}
\subsubsection{Optimal IP/IQP solutions}
\label{sec:mip_solutions}
We compute the optimal values $\Zopt$ of the OVRP with its IP formulation in \eqref{totalcost} and the ASP with its IQP formulation in \eqref{ASP-1} via the standard solver \textsc{CPLEX}, see \cref{sec:ex_setup}.
The optimal values solve as benchmark for the SDP relaxations of the QUBO formulations.
In \cref{tab:MIP_results} we provide the optimal values along with the corresponding running times.
The ASP problems are solved to optimality within a few seconds.
The \texttt{OVRP11} problem is solved to optimality within about two minutes, and \texttt{OVRP15} under 20 minutes.

	\begin{table}[h!]
		\begin{center}
			\begin{tabular}{|l||c|c|c||c|c|}
            \hline
				Instance &	 \texttt{ASP30} &  \texttt{ASP60} &  \texttt{ASP90} &    \texttt{OVRP11} &    \texttt{OVRP15} \\
                \hline
                \textbf{Optimal value $\Zopt$}  &	\textbf{246.78} & \textbf{820.08} & \textbf{1373.15} & \textbf{6254.19} & \textbf{3036.54} \\\hline
                Running time \textsc{CPLEX} in s &$0.45$ &$2.07$ &$2.82$ &$123.81$ &$1128.14$\\
				 \hline
			\end{tabular}
		\end{center}
        \caption{Running times and optimal values $\Zopt$ to IP and IQP formulations for the OVRP and ASP problems, respectively, using \textsc{CPLEX}. }
        \label{tab:MIP_results}
	\end{table}
The standard methods for IPs/IQPs implemented in solvers such as \textsc{CPLEX} are highly optimized for solving exactly this type of problem formulation \cite{junger2010integer}, yielding fast running times for these problems, see \cref{tab:MIP_results}. 
This is the result of decades of development and refinement, specifically tailored to these types of problems \cite{KOCH2022100031}.
Moreover, we emphasize here, that the corresponding problem instances were chosen such that their IP/IQP formulation is solvable to optimality providing a benchmark for comparisons. 
\subsubsection{QUBO Solver Results}
\label{sec:qubo_solver_results}
The QUBO reformulation, see \cref{sec:QUBO_reformulation}, changes the problem type from IP/IQP to QUBO, even though the optimal value $\Zopt$
 is the same.
 In the upcoming sections we compare different SDP relaxations of these QUBO formulations and use the optimal value $\Zopt$ as a benchmark. 
 Thus, for completeness, in the remainder of this subsection, we try to solve the QUBO formulations to optimality. 
 We point out that the QUBO reformulation is more difficult to solve than the original IP/IQP formulations due to the increased number of auxiliary binary variables, which in turn expand the search space of the problem, see \cite{gabbassov2024lagrangiandualityquantumoptimization}. Moreover, the penalty term in the QUBO formulation plays a significant role in the optimization process, see \cref{sec:binarization}, \cite{10.1007/s10878-014-9734-0} and \cite{alessandroni2024alleviatingquantumbigmproblem}. 
 
For reference and comparison with our goal of benchmarking the SDP relaxations, we provide the solutions to the QUBO formulation that we obtain from quadratic solvers implemented in
\textsc{CPLEX} and \textsc{GUROBI}.
There exists more specialized academic software such as \textsc{BiqMac} \cite{BiqMac}, \textsc{BiqCrunch} \cite{BiqCrunch}, \textsc{MADAM} \cite{MADAM}, \textsc{BigBin} \cite{BiqBin} as well as \textsc{QuBowl} \cite{QuBowl} and \textsc{McSparse} \cite{McSparse} for large sparse problems, specifically designed for solving QUBOs.
For the scope of this paper, we chose the quadratic solver included in \textsc{CPLEX} and \textsc{GUROBI} to serve as benchmark of QUBO solvers and direct the reader to the above articles for detailed comparisons of the QUBO solvers. 
Some of the solvers are only available as web services and on many benchmark instances used in \cite{QuBowl} the commercial solver \textsc{GUROBI} outperforms them.

Whereas for the IP/IQP formulations the solver stops with proven optimality, for the QUBO formulation we have to stop the solving process at a certain imposed time limit. Then we verify whether the constraints are satisfied by substituting the solution into the original problem and assess its optimality by comparing the objective value of the QUBO to $\Zopt$. 

\begin{figure}[h!]
    \centering
    \includegraphics[width=0.6\linewidth]{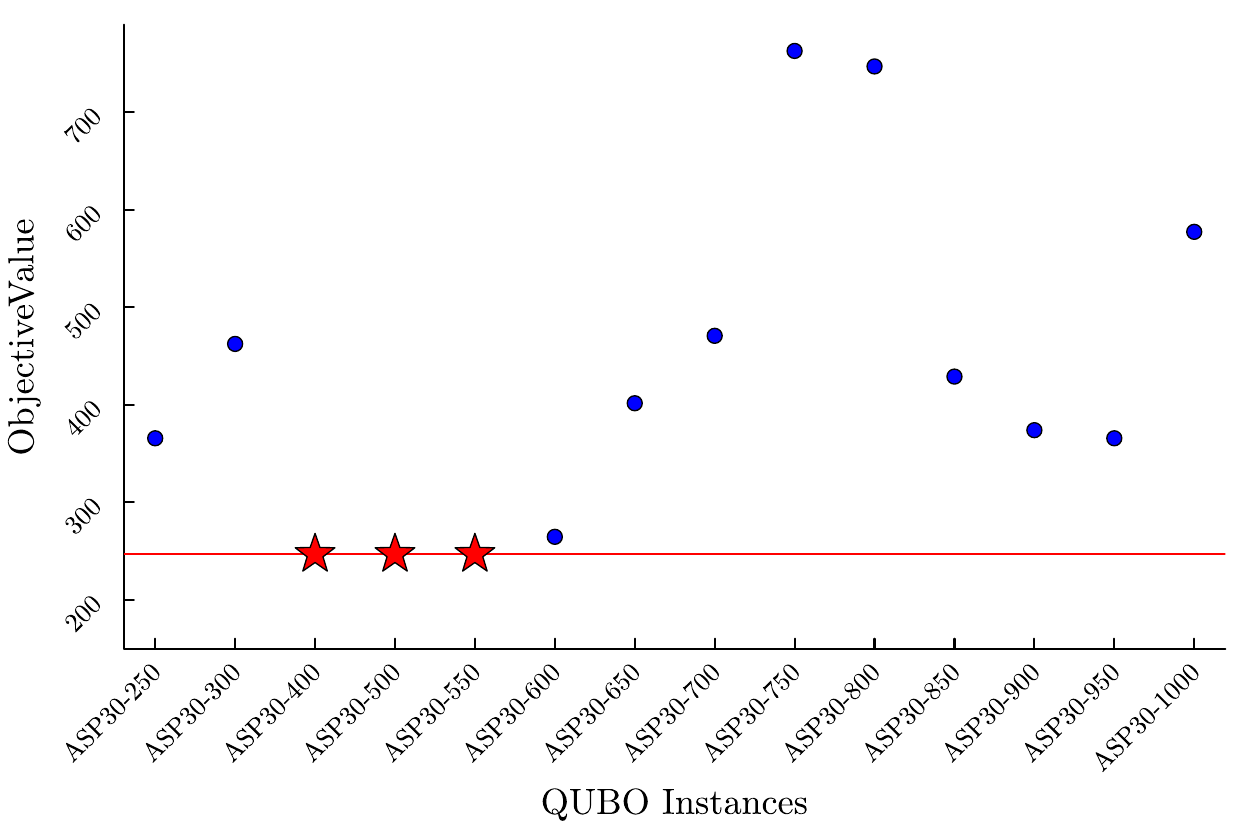}
    \caption{\textbf{Objective values} $\zQUBO$ for solving the \textbf{QUBO formulations} of \texttt{ASP30} with \textbf{different penalties} using \textsc{CPLEX} and a \textbf{time limit of 120 seconds}. The red line indicates the optimal value $\Zopt$ of \texttt{ASP30}. The blue dots are the objective values $\zQUBO$, additionally the red stars indicate that the optimal value was reached for these formulations.}
    \label{fig:QUBO_ASP30_120s}
\end{figure}

We denote the objective value reached by solving the QUBO formulation after a certain time limit as $\struc{p_{\text{QUBO}}}$. 
This is an \textbf{upper bound} to $\Zopt$.
In \cref{fig:QUBO_ASP30_120s}, we exemplarily display $\zQUBO$ that we obtain for the  \texttt{ASP30} instances with a time limit of 120 s with the solver \textsc{CPLEX}. 
For the penalties 400, 500 and 550, optimality is reached.
Note that this is an upper bound for the time needed to reach the optimal value.
For the other instances we find feasible solutions, but not the optimal one, after 120 s. 
We see here that the chosen penalty in the QUBO construction influences the result significantly.

We denote the absolute difference of the upper bound $\zQUBO$ to the optimal value $\Zopt$ and by the relative difference, respectively, by:
\begin{align*}
\struc{\DeltaQUBOabs} =|\zQUBO-\Zopt| \text{ and } \struc{\DeltaQUBO} = \frac{\DeltaQUBOabs}{\Zopt}
\end{align*} 
For  \texttt{ASP60} and  \texttt{ASP90}, we find feasible solutions after a set time limit of $120 s$ using \textsc{CPLEX}.
However, these upper bounds are very large, e.g.~the relative difference to the optimal value is $\DeltaQUBO>1300\%$ for \texttt{ASP90}-1450, meaning that the upper bound is by a factor of 13 larger than the optimal value.
For the OVRP instances, we did not obtain a feasible solution after setting the time limit to 30 minutes using \textsc{CPLEX}. 
This indicates the difficulty of the QUBO formulation compared to the IP/IQP formulations that we could solve to optimality, see \cref{tab:MIP_results}.

We additionally try to tackle the problems via the standard solver \textsc{GUROBI} \cite{gurobi} that produces better solutions for the larger instances  \texttt{ASP60} and  \texttt{ASP90}.
We give an exemplary overview on $\DeltaQUBO$ for different time limits in \cref{tab:QUBO_example_optimality}, see \cref{tab:QUBO_example_optimality_abs_difference} for $\DeltaQUBOabs$ correspondingly.

\begin{table}[h!]
\begin{center}
    \begin{tabular}{|l||c|c||c|c||c|c|}
    \hline
        Instance & Time Limit & $\DeltaQUBO$ & Time Limit & $\DeltaQUBO$ & Time Limit & $\DeltaQUBO$ \\
        \hline
         \texttt{ASP30}-$950$ & $10s$ & $0$ &$120s$ &$0$ &$600s$ &$0$\\
         \texttt{ASP60}-$900$ & $10s$ & $0.29$ &$120s$ & $0.29$&$600s$ &$0.29$\\
         \texttt{ASP90}-$1450$ & $10s$ & $0.40$ &$120s$ &$0.21$ &$600s$ &$0.19$\\
        \hline
    \end{tabular}
\end{center}
\caption{\textbf{Relative difference} $\DeltaQUBO$ of the upper bound $\zQUBO$ to the optimal value for selected \textbf{QUBO}s using \textsc{GUROBI} with different time limits. For \texttt{ASP30}-950 we find the optimal value. The provided time limits are upper bounds on the running time to reach the respective objective value $\zQUBO$. The computation is stopped at the provided time limit.}
\label{tab:QUBO_example_optimality}
\end{table}

In particular, for the instance  \texttt{ASP60}-900, we have $\DeltaQUBO\approx 29\%$.
Comparing to \cref{fig:QUBO_ASP30_120s}, different penalties work well on \textsc{GUROBI}. 
While we only find a feasible, but not optimal solution for \texttt{ASP30} with penalties $\{400, 500, 550\}$, we find an optimal solution for  \texttt{ASP30}-950.
 
Similarly to \textsc{CPLEX}, \textsc{GUROBI} does not provide feasible solutions for either OVRP instance within a time limit of 30 minutes.
For reference, in \cref{sec:dwave_results}, we exemplarily provide results using the D-Wave's Quantum Annealer Procedure that we access via \textsc{GAMS}'s \textsc{QUBO\_SOLVE} tool, see \cref{sec:ex_setup}.
Similarly, we find the optimal solution for \texttt{ASP30}, and feasible solutions for \texttt{ASP60} and \texttt{ASP90} but the OVRP instances give infeasible solutions.

\subsection{Experimental Results of SDP relaxations}
\label{sec:results}
In the following \cref{sec:results_vrp,sec:results_asp} we provide the results of the SDP relaxations concerning the larger QUBOs emerging from the OVRPs and from the smaller QUBOs from the ASPs, respectively. 
We tackle these QUBOs via the SOS-SDP relaxations from \cref{sec:SOSSDP} and HU from \cref{sec:HU}.

\subsubsection{SDP relaxations of QUBO formulations for OVRP}
\label{sec:results_vrp}
In \cref{tab:VRPresults} we summarize the results of the SDP relaxations for the QUBOs constructed from the OVRPs for SOS-SDP and HU methods. 
	\begin{table}[h]
		\begin{center}
			\begin{tabular}{|l|c|c|c|c|c|c|}\hline
				&	\multicolumn{3}{c|}{\textbf{\texttt{OVRP11}}} & \multicolumn{3}{c|}{\textbf{\texttt{OVRP15}}} \\
				
				&\multirow{2}*{} \multirow{2}*{$\DeltaSOS|\DeltaHU$} & $ \tSOS|\tHU$ & $\tHUQuantum$ & \multirow{2}*{$\DeltaSOS|\DeltaHU$ }&$ \tSOS|\tHU$ & $\tHUQuantum$  \\
                & & in s & in s & & in s & in s\\\hline 
                \rowcolor{gray!20} \textsc{SOS.jl}   & $\dagger$& $\dagger$& - &  $\dagger$& $\dagger$ & - \\
				\rowcolor{gray!20} \textsc{TSSOS}  & $\dagger$ & $\dagger$& - & $\dagger$& $\dagger$ & - \\
				\rowcolor{gray!20} \textsc{ManiSDP}   & $1.32$ & $28.08$& - & $1.60$ & $95.63$ & -\\
				HU\tiny{$\epsilon = \num{e-2}$} &$ \num{7.9e6}$ & $58.01$ & $ \num{9.3e12}$ & $ \num{5.6e7}$ & $310.36$
					& $ \num{9.4e12}$ \\
				HU\tiny{$\epsilon = \num{e-3}$} &$ \num{7.4e5}$ & $348.02$ & $\num{3.3e18}$ & $\num{6.2e6}$ & $1314.95$ & $\num{2.0e18}$\\\hline        
			\end{tabular}
		\end{center}
        \caption{\textbf{Results of SDP relaxations} of QUBOs referring to \texttt{OVRP11} and \texttt{OVRP15} comparing SOS and HU methods. Rows for SOS methods are gray, rows for HU methods are white. We display relative differences of the lower bounds to the optimal value by $\DeltaSOS$ and $\DeltaHU$, respectively, as well as running times for solving the SDP by $\tSOS$ and $\tHU$, respectively, depending on the method, indicated by '$|$'. Note that $\tHUQuantum$ is not available for SOS methods, as indicated by '-'. A dagger '$\dagger$' indicates that no result is available due to solver break down or exceeding maximal running time. }
        \label{tab:VRPresults}
	\end{table}

\textbf{SOS-SDP.} From both \textsc{SOS.jl} and \textsc{TSSOS} we do not receive a result for any of the QUBOs resulting from   \texttt{OVRP11} and \texttt{OVRP15}.
We stop \textsc{COSMO} inside \textsc{TSSOS} after four hours for OVRP11. 
All internal SDP solvers, i.e., \textsc{Mosek} and \textsc{Loraine}, as well as \textsc{SCS} break down for the first order relaxation.

Using ManiSDP we are able to execute the first order relaxation for both problems and produce a lower bound within approximately 0.5 and 1.5 minutes, that still yield a relative difference $\DeltaSOS$ of more than 100\% of the optimal IP values $\Zopt$.
These results from \textsc{ManiSDP} are the only benchmarks we achieve on the SOS side, concerning the OVRP instances.

\textbf{HU.} We execute HU for the precision $\epsilon$ in $\{\num{e-2}, \num{e-3}\}$ and produce lower bounds $\zHU$ in a few minutes for \texttt{OVRP11}, and in about 20 minutes for \texttt{OVRP15}, using $\epsilon=\num{e-3}$. 
Translated to the quantum version of HU, based on the assumptions in \cref{sec:experiments_attributes} and \cite{henze2025} it would take more than $10^{10}$ years to solve the SDP on a quantum computer with $\epsilon=\num{e-3}$ for both problems.
Even then, the lower bounds achieved are extremely far away from the optimal value.

The above results for HU are performed on a CPU. 
Using a GPU, we can accelerate the running time $\tHU$ of the classical HU version such that we can use smaller tolerances $\epsilon$ in a reasonable time frame. 
We present the corresponding results in \cref{tab:VRPresults_GPU,tab:VRPresults_GPU_upperbounds}.

\begin{table}
		\begin{center}
			\begin{tabular}{|l|c|c|c|c|c|c|}\hline
				&	\multicolumn{3}{c|}{\texttt{OVRP11}} & \multicolumn{3}{c|}{\texttt{OVRP15}} \\
				
				& \multirow{2}*{$\DeltaHU$} & $\tHU$& $\tHUQuantum$ & \multirow{2}*{$\DeltaHU$} &$\tHU$& $\tHUQuantum$ \\
                & & in s & in s & & in s & in s \\
				\hline
				HU$^*$\tiny{$\epsilon = \num{e-2}$} &$ \num{7.9e6}$ & $\mathbf{12.51}$ & $ \num{9.3e12}$ & $ \num{5.6e7}$ & $\mathbf{72.64}$
					& $ \num{9.4e12}$ \\
				HU$^*$\tiny{$\epsilon = \num{e-3}$} &$ \num{7.4e5}$ & $\mathbf{68.83}$ & $\num{3.3e18}$ & $\num{6.2e6}$ & $\mathbf{266.42}$ & $\num{2.0e18}$\\
                HU$^*$\tiny{$\epsilon = \num{e-4}$} &$ \num{7.1e4}$ & $\mathbf{16578.44}$ & $\num{2.3e24}$ & $\num{5.5e5}$ & $\mathbf{37405.86}$ & $\num{4.7e24}$\\
                \hline
			\end{tabular}
		\end{center}
        \caption{\textbf{Results of SDP relaxations via HU} of QUBOs referring to \texttt{OVRP11} and \texttt{OVRP15}. In contrast to \cref{tab:VRPresults} these results are obtained via execution on a GPU, indicated by '$*$' for distinction. Running times $\tHU$ are marked bold.}
        \label{tab:VRPresults_GPU}
	\end{table}

We see that the GPU decreases running times significantly.
However, solving \texttt{OVRP15} with $\epsilon=\num{e-4}$ takes approximately 10 h with a relative difference $\DeltaHU\approx \num{5.5e5}$ of the lower bound $\zHU$ to the optimal value.
Through the randomized rounding procedure in HU, we also produce upper bounds $\zHUupper$ that are closer to the optimal value, e.g.~we have $\DeltaHUupper\approx\num{2.2e4}$ for \texttt{OVRP15} at $\epsilon = \num{e-4}$.
Nevertheless, this is still very far from the optimal value.

We summarize that the QUBOs derived from the OVRPs are too large to be tackled via the SDP solvers within the SOS methods. While we can at least solve the problem with \textsc{ManiSDP} within 0.5 to 1.5 min, the second order relaxation of the Lasserre hierarchy cannot be executed, leaving still a large relative difference $\DeltaSOS$ to the optimal value.
The HU method can tackle the large problem producing lower and upper bounds but for reasonable running times only with low precision.

\subsubsection{SDP relaxations of QUBO formulations for ASP}
\label{sec:results_asp}
Compared to the OVRP instances, the ASP instances can be tackled by the SOS-SDP methods. 
We first describe the results for the SOS-SDP methods and afterwards for HU.

\textbf{SOS-SDP.} 
	\begin{table}[h]
	\begin{center}
		\begin{tabular}{|l|c|c|c|c|c|c|}\hline
			&	\multicolumn{2}{c|}{\texttt{ASP30}} & \multicolumn{2}{c|}{\texttt{ASP60}}  & \multicolumn{2}{c|}{\texttt{ASP90}} \\
			&\multirow{2}*{} \multirow{2}*{Average} & \multirow{2}*{Average} & \multirow{2}*{Average} & \multirow{2}*{Average}& \multirow{2}*{Average}& \multirow{2}*{Average}  \\[+4pt] 
			& $\DeltaSOS$ & $\tSOS$ &$\DeltaSOS$ &$\tSOS$ &$\DeltaSOS$ & $\tSOS$ \\
			\hline
			\textsc{ManiSDP} & $1.99$& $1.45 s$& $3.01$ &  $0.82 s$& $2.41$ & $1.45 s$\\\hline 
			\textsc{SOS.jl} w. \textsc{Mosek}  & $1.99$& $9.91 s$& $3.01 $ &  $221.25 s$& $2.41$ & $8256.9 s$ \\\hline   
			\textsc{SOS.jl} w. \textsc{Loraine}  & $1.99$& $670.66 s$& $\dagger$ &  $\dagger$& $\dagger$ & $\dagger$\\\hline 
			\textsc{TSSOS} w. \textsc{Mosek} & $1.99$& $3.24 s$& $3.01$ &  $112.42 s$& $\dagger$ & $\dagger$\\\hline 
			\rowcolor{gray!20}\textsc{TSSOS} w. \textsc{Mosek} & $0.13$& $28.46 s$& $0.21$ &  $912.47s$& $\dagger$ & $\dagger$\\\hline 
		\end{tabular}
	\end{center}
	\caption{\textbf{Results of SOS-SDP relaxations} of QUBOs referring to \texttt{ASP30}, \texttt{ASP60} and \texttt{ASP90}. We display the \textbf{average} of the relative difference $\DeltaSOS$ to the optimal value and running times $\tSOS$ for each set of instances. The gray row refers to the second order Lasserre relaxation, the other rows refer to the first order relaxation. A dagger '$\dagger$' indicates that no result is available due to solver break down or exceeding maximal running time. For \textsc{Mosek}, we display the result with the default tolerance tol=$10^{-8}$. For \textsc{SOS.jl} with \textsc{Mosek} the total solving time of all \texttt{ASP90} instances exceeds our set time limit. We display a rerun for these instances without time limit.}
	\label{tab:ASP_running_times}
\end{table}
We give an overview on the most promising results in \cref{tab:ASP_running_times}, see also \cref{tab:ASP_running_times_abs} in the appendix, by presenting the average over the relative differences $\DeltaSOS$ and running times $\tSOS$ over each set \texttt{ASP30},  \texttt{ASP60} and \texttt{ASP90} with different penalties, respectively.

    \begin{figure}[h]
    \centering
    \includegraphics[width=0.7\linewidth]{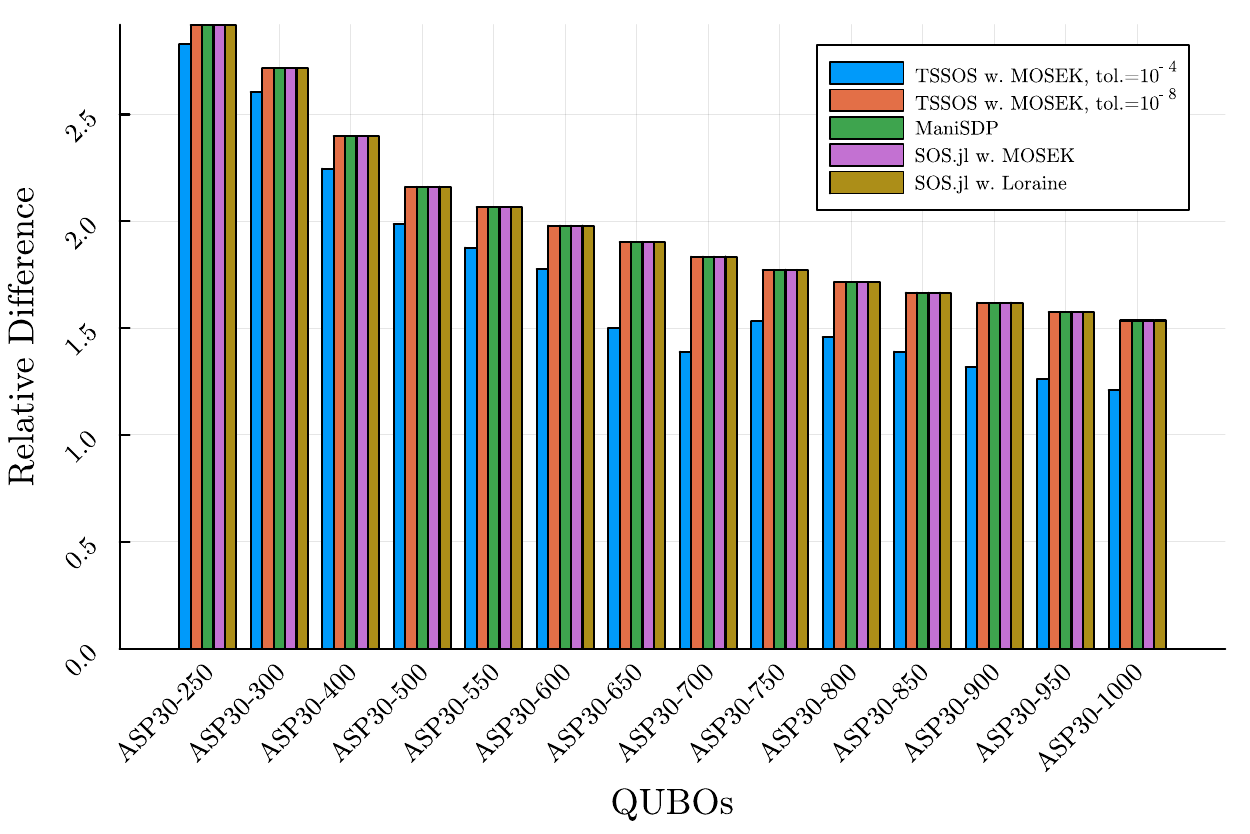}
    \caption{{Relative Difference} $\DeltaSOS$ on all \texttt{ASP30} instances with different penalties, obtained with \textbf{first order} relaxations, $2d=2$. 
    }
    \label{fig:ASP30_relax2}
\end{figure}
    
In particular, we are able to receive solutions for a range of methods on the smallest instance  \texttt{ASP30}. 
In \cref{fig:ASP30_relax2}, we show the relative difference $\DeltaSOS$ to the optimal value obtained by different solvers for each penalty. 
This display resembles the first order Lasserre relaxation with a relative difference to the optimal value of $\DeltaSOS > 100\%$ for almost all solvers and instances.
Interestingly, the lower bounds improve slightly with a a higher penalty in the QUBO formulation.
We observe that using the default parameters in \textsc{Mosek}, \textsc{Loraine} and \textsc{ManiSDP}, we get the same lower bounds. 
We note that, while \textsc{Mosek} reports a FEASIBLE\_POINT, \textsc{Loraine} reports an UNKNOWN\_RESULT\_STATUS for these bounds.
The fastest solver is ManiSDP with an average solving time of about $1.5 s$ for all  \texttt{ASP30} instances.
The solver \textsc{COSMO} inside \textsc{TSSOS} does \textit{not} output a lower bound for any instance, the \textsc{SCS} solver on the Ising formulation does not reliably give lower bounds for all penalties, using the default settings. 

\begin{figure}[h!]
    \centering
     \includegraphics[width=0.75\linewidth]{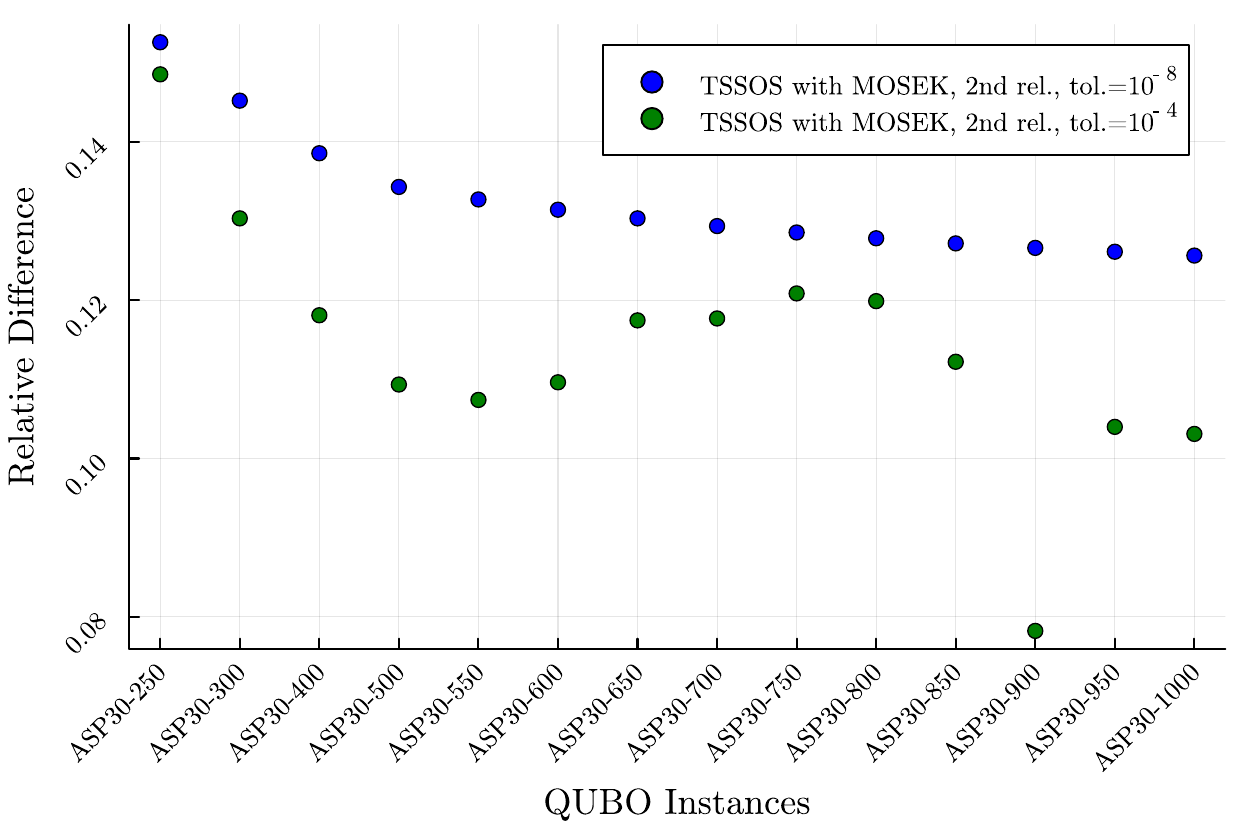}
    \caption{{Relative Difference} $\DeltaSOS$ on all \texttt{ASP30} instances with different penalties, obtained with \textbf{\textsc{TSSOS}}, using \textsc{Mosek} with the \textbf{second order} relaxation, $2d=4$ and the default tolerances of $\num{e-8}$, as well as $\num{e-4}$.}
        \label{fig:ASP30_MOSEK2ndTSSOS}
\end{figure}

For \texttt{ASP30}, we are able to execute the second order Lasserre relaxation via \textsc{TSSOS} and the solver \textsc{Mosek}\footnote{Both \textsc{ManiSDP} and \textsc{TSSOS} with \textsc{COSMO} can also execute the second order relaxation, however the results are not promising: 
\textsc{COSMO} does not yield lower bounds and \textsc{ManiSDP} reaches lower bounds with $\DeltaSOS$ of order $10^5\%$ in an hour of running time. Using \textsc{Mosek} and \textsc{COSMO} inside \textsc{SOS.jl}, the solver breaks down for the second order relaxation.}. 
We display the results for two different parameter settings in \cref{fig:ASP30_MOSEK2ndTSSOS}. 
The second order relaxation provides close lower bounds to the optimal value:
On average, we can achieve a relative difference $\DeltaSOS\approx 13\%$ ($\DeltaSOSabs\approx 32.7$) using \textsc{TSSOS} with \textsc{Mosek} and $\text{tol}~=~10^{-8}$ with a running time of $28.5$s. 
We do not achieve a higher relaxation order.
Note that we see from \cref{fig:ASP30_MOSEK2ndTSSOS} that we get slightly better bounds with a larger tolerance of $\text{tol}~=~10^{-4}$, though these bounds are not certified as \textsc{FEASIBLE\_POINT} in the solver output.

\begin{figure}[h]
	\centering
	\includegraphics[width=0.8\linewidth]{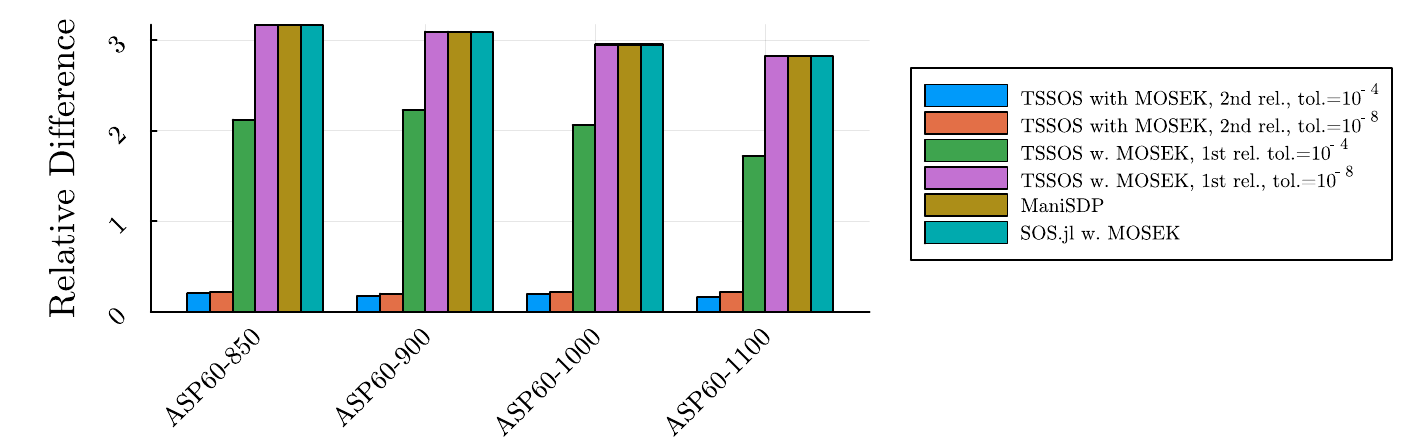}
	\caption{{Relative Difference} $\DeltaSOS$ on \texttt{ASP60} instances with different penalties. 
	}
	\label{fig:ASP_60_SOS}
\end{figure}

We give an overview on all \texttt{ASP60} instances in \cref{fig:ASP_60_SOS}.
We again receive the best lower bounds by running the second order relaxation with \textsc{TSSOS} using \textsc{Mosek}.
With the default tolerance, tol $=10^{-8}$, this yields an average relative difference $\DeltaSOS$ of $21\%$ ($\DeltaSOSabs\approx 174.3$) with a running time of over $15$ minutes on average. 

However, for the larger \texttt{ASP90} instances, we do not achieve any results with \textsc{TSSOS}, the internal solver \textsc{Mosek} breaks down even for the first order relaxation.
We can only tackle \texttt{ASP90} via the first order relaxation inside \textsc{SOS.jl} using  \textsc{Mosek} and \textsc{ManiSDP}. 
The former yields an average relative difference $\DeltaSOS$ of $240\%$ ($\DeltaSOSabs\approx 3308.5$) within an average running time $\tSOS$ of over $2$h while \textsc{ManiSDP} achieves the same result in on average $1.5$s.

\begin{figure}[h]
    \centering
    \includegraphics[width=0.8\linewidth]{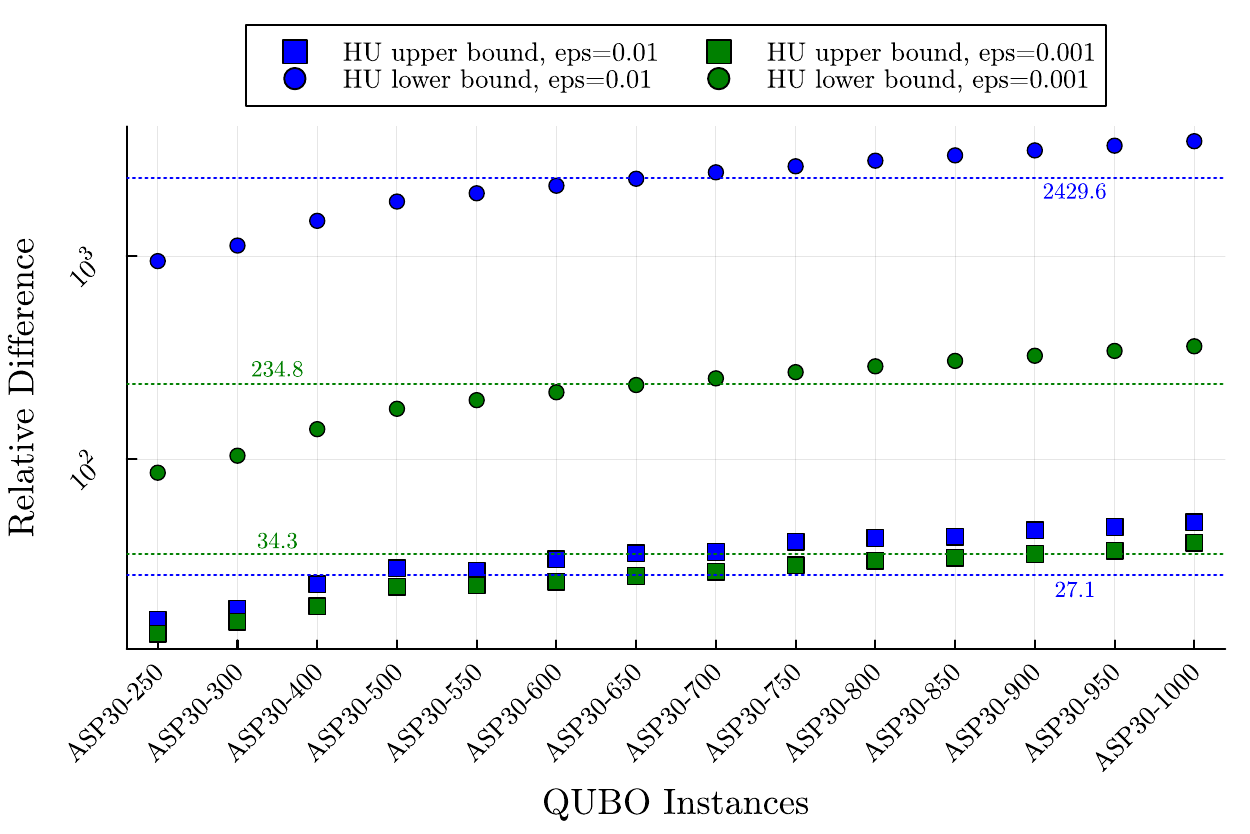}
    \caption{\textbf{Relative Differences} $\DeltaHU$, $\DeltaHUupper$ of \textbf{lower bounds} $\zHU$ and \textbf{upper bounds} $\zHUupper$, respectively, to optimal value, on all \texttt{ASP30} instances with different penalties, obtained with HU and precisions $\epsilon=\num{e-2}$ and $\epsilon=\num{e-3}$. Dotted lines display the average values for $\DeltaHU$ and $\DeltaHUupper$ over all \texttt{ASP30} instances. The $y$-scale is logarithmic.}
    \label{fig:ASP_HU_30}
\end{figure}

\textbf{HU.} We display the relative differences of the lower bounds as well as the upper bounds obtained by the randomized rounding procedure for all \texttt{ASP30} instances in \cref{fig:ASP_HU_30}.
We see that the upper bounds are closer to the optimal value than the lower bounds, yet still very far away. 
In particular, on average, the upper bounds to the optimal value yield $\DeltaHUupper\approx2700\%$ ($\DeltaHUabsupper~\approx~6688.7$) for precision $\epsilon=\num{e-3}$, while for the lower bounds, the average $\DeltaHU$ is of order $\num{e4}\%$.
We notice that for our given instances, the upper and lower bounds improve with a decreasing penalty in the QUBO formulation. 
This can be explained by the fact that the precision of the HU solution gets scaled with the norm of the cost matrix $C$ of the Ising problem \eqref{eq:ising}. As larger penalties lead to a higher norm of $C$, consequently the gap between lower and upper bound becomes larger.
Running the experiments on the default CPU, we produce these results within few seconds for $\epsilon = \num{e-2}$ and on average $57$ s per instance for $\epsilon=\num{e-3}$.

In \cref{tab:ASPresults_GPU_example,tab:ASP60_90results_GPU_example} we exemplarily compute the SDP relaxations via HU on a GPU for the instances \texttt{ASP30}-250, \texttt{ASP60}-850 and \texttt{ASP90}-1400.
For this exemplary display, we choose the instances with the smallest penalty from the QUBO construction.
The GPU provides faster results than the default CPU, letting us try smaller tolerances $\epsilon$ in a reasonable running time.
With $\epsilon=10^{-5.5}$, we receive a lower bound $\zHU$ for \texttt{ASP30}-250 that has a relative difference $\DeltaHU\approx 310\%$, taking a classical running time of more than three days. 
To achieve this precision for \texttt{ASP30}-250 using the quantum version of HU would take an order of $10^{20}$ years.

	\begin{table}[h]
		\begin{center}
			\begin{tabular}{|l|c|c|c|c|}\hline
				&	\multicolumn{4}{c|}{\textbf{\texttt{ASP30}-250}} \\
				
				& $\DeltaHU$ & $\DeltaHUupper$ & $\tHU$ in s & $\tHUQuantum$ in s  \\
                \hline
                HU$^*$\tiny{$\epsilon = \num{e-2}$} & $945.6$ & $16.8$ &$0.17$ & $\num{8.4e10}$\\
                 HU$^*$\tiny{$\epsilon = \num{e-3}$} &$86.0$ & $14.2$ &$3.08$ &$\num{1.2e16}$ \\
                  HU$^*$\tiny{$\epsilon = \num{e-4}$} &$8.8$ & $6.6$ &$500.36$ &$\num{3.1e21}$ \\
                   HU$^*$\tiny{$\epsilon = \num{e-5}$} &$3.8$ &$3.2$&$21861.1$ &$\num{7.7e25}$ \\
                    HU$^*$\tiny{$\epsilon = 10^{-5.5}$} &$3.1$ & $2.8$&$\num{2.7e5}$ &$\num{3.5e28}$ \\
                    \hline
			\end{tabular}
		\end{center}
        \caption{\textbf{Results of SDP relaxations via HU} of QUBOs referring to \texttt{ASP30-250}. These results are obtained via execution on a \textbf{GPU}, indicated by '$*$' for distinction. We display relative differences to optimum $\DeltaHU$ and $\DeltaHUupper$ of lower bounds and upper bounds, respectively.}
        \label{tab:ASPresults_GPU_example}
	\end{table}

    	\begin{table}[h]
		\begin{center}
			\begin{tabular}{|l|c|c|c|c|c|c|}\hline
				&	\multicolumn{3}{c|}{\texttt{ASP60-850}} &	\multicolumn{3}{c|}{\texttt{ASP90-1400}} \\
				
				& $\DeltaHU$ & $\tHU$ in s & $\tHUQuantum$ in s & $\DeltaHU$ & $\tHU$ in s & $\tHUQuantum$ in s  \\
                \hline
                HU$^*$\tiny{$\epsilon = \num{e-2}$} & $\num{5.4e3}$ &$0.23$ & $\num{2.6e11}$& $\num{1.9e4}$ &$0.5$ & $\num{5.2e11}$ \\
                 HU$^*$\tiny{$\epsilon = \num{e-3}$} &$807.5$ &$3.98$ &$\num{4.2e16}$ & $\num{2.7e3}$ &$4.99$ & $\num{6.7e16}$\\
                   HU$^*$\tiny{$\epsilon = \num{e-4}$} &$74.42$ &$504.7$ &$\num{1.1e22}$ & $\num{256.2}$ &$1313.8$ & $\num{5.8e22}$\\
                    \hline
			\end{tabular}
		\end{center}
        \caption{\textbf{Results of SDP relaxations via HU} of QUBOs referring to \texttt{ASP60-850} and \texttt{ASP90-1400}. These results were obtained via execution on a \textbf{GPU}, indicated by '$*$' for distinction. See \cref{tab:ASP60_90results_GPU_example_upper} for upper bounds.}
        \label{tab:ASP60_90results_GPU_example}
	\end{table}

\begin{figure}[h]
    \centering
    \includegraphics[width=0.75\linewidth]{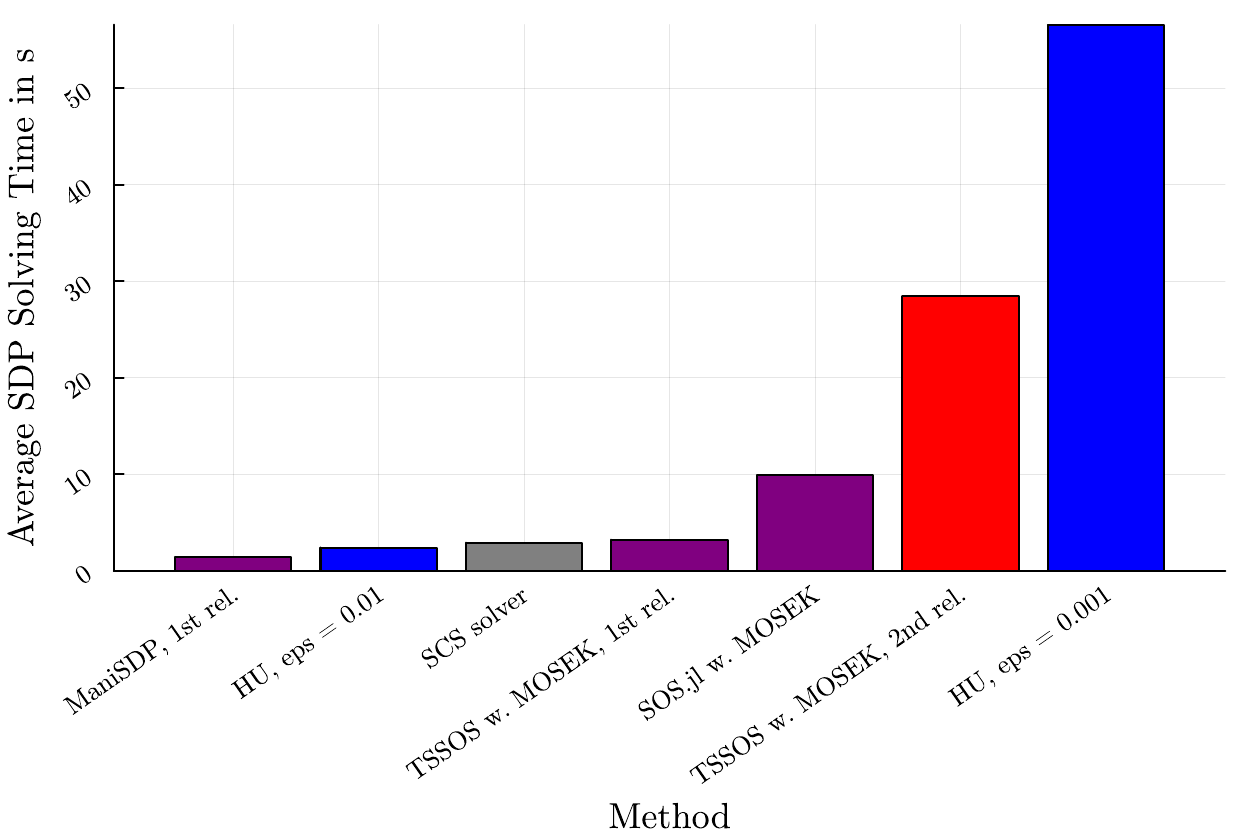}
    \caption{\textbf{Averages} over \textbf{running times} $\tSOS$ and $\tHU$ of all \texttt{ASP30} instances for different methods. 
    The red bar indicates that the average relative difference $\DeltaSOS\leq 0.2$. The purple bars indicate that $\DeltaSOS\leq 2$. The blue bars indicate that $\DeltaHU$ is larger than 2. The gray bar indicates, that the solver did not output a reliable lower bound $\DeltaSOS$. For \textsc{TSSOS} with \textsc{Mosek} tol=$10^{-8}$ is displayed. If there is no other explicit description, the first order Lasserre relaxation is displayed for SOS methods. \textsc{SOS.jl} with \textsc{Loraine} reached an average relative difference $\DeltaSOS\approx 1.99$ but is not displayed due to average running time of $>600 s$. HU is executed on the default CPU.}
    \label{fig:ASP_running_times}
\end{figure}

We summarize our findings on both SOS and HU methods, displaying the \texttt{ASP30} instances in \cref{fig:ASP_running_times}: The plot shows the average running times $\tSOS$ and $\tHU$ for several methods, while also indicating the relative differences $\DeltaSOS$ and $\DeltaHU$. 
Only via the second order Lasserre hierarchy with \textsc{TSSOS} using \textsc{Mosek}, we could achieve $\DeltaSOS<20\%$. 
The fastest running times are performed by \textsc{ManiSDP} and the classical version of HU with low precision.

\subsection{Discussion}
\label{sec:discussion}
We highlight again that the goal of this study is not to suggest or promote a certain method but to compare existing methods for SDP relaxations on \textbf{real-world data}.

The parameter settings in the real-world applications are chosen with the intention to make the instances solvable for IP/IQP solvers and thus provide optimal benchmarks.
As described in \cref{sec:industryproblems}, the downside is that these settings yield the smallest instances to be encountered in industry.
However, from the -- for industry standards -- small instances, we observe that the \textbf{QUBO formulations} of these problems are hard for all tested methods. 

As expected, the standard solvers \textsc{GUROBI} and \textsc{CPLEX} optimally solve the IQP respectively IP formulations of the smaller ASP instances in seconds and the largest OVRP instance in under 20 minutes.
As shown in the QUBO benchmarks by Mittelmann \cite{mittelmann}, specialized QUBO solvers, see \cref{sec:qubo_solver_results}, can solve the QUBO instances from the QPLIB \cite{qplib} library with up to 700 variables, and sparse instances even up to 10.000 variables \cite{QuBowl}.
Especially, dense QUBO instances are often randomly generated \cite{wiegele2007biq,QuBowl,wang2024solvinglowranksemidefiniteprograms,loraine2023}.
With the ASP instances, we provide QUBO formulations from real-world instances with about 50\% non-zero entries.
Recall, that the largest ASP instance \texttt{ASP90} has 276 variables.
Hence, merely based on the number of variables, it is expectable that the sizes of our instances are solvable for specialized QUBO solvers.
We point out here, that our instances are not QUBO instances by nature nor MaxCut instances, and we reformulate their IP/IQP formulations to QUBOs.
From \cref{sec:qubo_solver_results}, we see that, using standard solvers \textsc{CPLEX} and \textsc{GUROBI}, we can only solve the smaller instance \texttt{ASP30} with 96 variables to optimality for certain penalties in the QUBO reformulation, while we obtain no feasible solution for the OVRP instances.
While we do not use the specialized solvers in the scope of this paper, we see  the difficulty of the reformulated problems from the behavior of the standard solvers.

The focus on our paper lies on \cref{sec:results}, where we relax the QUBO formulations and provide lower bounds to the optimal value.
Using the \textbf{SOS-SDP} approach, we test several software from polynomial optimization.
As expected, we observe that the SDP solvers in their default settings do not always return reliable solutions, see e.g.~\cite{Sremac}, do not output lower bounds to the optimal values, or even break down given the size of the SDPs.
The software \textsc{TSSOS} with internal SDP solver \textsc{Mosek} provides the qualitatively best lower bound for the smallest problems \texttt{ASP30} and \texttt{ASP60} with an relative difference to the optimal value of 13\% and 21\% respectively, performing the second Lasserre relaxations. 
However, we cannot tackle the larger problems and higher relaxations.
Using the other SOS-SDP methods, we can only provide the first order relaxations that give lower bounds with relative differences to the optimal value of $\gg100\%$ for the tested instances.
In terms of running time, the SDP solver \textsc{ManiSDP} is very fast and finds these bounds in very few seconds for the ASP instances, as well as being the only tested SOS-SDP method that gives results for the OVRP instances, within only 1.5 minutes.

\textbf{Hamiltonian Updates} is an algorithm that scales in favor of the matrix input size, providing the opportunity to solve very large SDPs. 
However, this comes with a costly dependence of the running time with the precision. 
As expected, we can solve all instances via HU within seconds or few minutes for a low precision.
However, on the tested instances, we can only give very poor lower bounds to the optimal value, where the best achieved relative differences to the optimal values is still $\gg 300\%$.

The running time of the quantum version of HU can be approximated via counting the quantum gates in an optimistic way, in favor of the physical realization of a quantum computer. 
We see, however, that the approximated running times exceed its classical version tremendously on our tested instances, see \cite{henze2025} for a non-asymptotic analysis of the running time, showing that the quantum approach does not beat the classical version on realistic instances despite asymptotical advantage.

Using the randomized rounding procedure, we furthermore get upper bounds to the optimal value, as well as a feasible solution vectors. 
While these upper bounds are also far away from the optimal values, it is an advantage over most current implementations of SOS-SDP methods, where solutions extraction is based on the Fialkow-Curto approach, which is only guaranteed to be applicable if the hierarchy converges to the optimal value; see e.g.~\cite{Laurent:Survey}.

\section{Conclusion}
\label{sec:conclusion}
It is widely known that in optimization real-world data behaves differently than randomly generated data, which is often used in the development of algorithms and software.
In this work, we make the same observation. 
We see moreover that reformulating IP/IQP problems to QUBOs increases the difficulty of the instances if they are not naturally QUBO or MaxCut instances.
Our test instances are by far best solved by the specialized and well-developed IP/IQP methods.

The goal of this work was to compare SDP relaxations of QUBO formulations stemming from real-world IPs/IQPs.
We see that we can solve the smallest QUBOs directly, using quadratic solvers of \textsc{GUROBI} and \textsc{CPLEX}, and still achieve better lower bounds than with the first order SOS-SDP relaxations and the SDP relaxation via HU.
If only comparing the SDP approaches, we observe that on the chosen instances the quality of the lower bounds provided by SOS-SDP is better than using HU. 

While on the tested instances the quality of solutions is non-competitive, HU can provide lower bounds to very large problems, where SOS-SDP methods and QUBO solvers fail, and it also gives upper bounds and feasible solutions.
Generally, the practitioner needs more than just a lower bound on the solution -- i.e., an actual solution to the optimization problem is needed, or at the very least, upper and lower bounds. 
Even though HU scales unfavorably with high precision, it does provide these properties. Therefore, it would be interesting to  conduct further research on its applicability to large real-world instances.

In general we notice that also on the polynomial optimization side, the results are very much dependent on the method and the solver.
Especially, approaches exploiting sparsity perform well in our experiments, directing the future research to more specialized solvers and software for preprocessing specifically considering the structure and properties of the problem.
With this work, we confirm the trend that developing these methods is fruitful, as well as testing them not only on randomly generated data.

\section{Acknowledgments}
\label{sec:ackknowledgements}
All authors were supported by the German Federal Ministry for Economic Affairs and
Climate Action (BMWK), project ProvideQ. 
FH and DG are also supported by the
German Federal Ministry of Education and Research (BMBF), project QuBRA and by Germany's
Excellence Strategy -- Cluster of Excellence Matter and Light for Quantum Computing (ML4Q) EXC 2004/1 (390534769). 
We thank Sabrina Ammann for helpful comments on the draft, as well as Jie Wang, Victor Magron, Amir Ali Ahmadi and Christoph Helmberg for insightful discussions.

\section{Authors' Contributions}
TG, WG, FF, TdW and DG initiated the project.
TdW lead the project and coordinated it jointly with BO.
TdW, DG, FF and WG supervised the project locally in the respective research groups.
BO, FH and VJ carried out the implementation of the different methods, ran preliminary experiments and conceptualized the experimental study, AD assisted them.
AD implemented unifying scripts of the codes to run the experiments. 
BO, FH, VJ evaluated the results and TdW, DG, FF and TG gave feedback on the results.
AD, BO, VJ and FH wrote the code documentation.
BO, TG, VJ, FH and TdW wrote the initial draft of the manuscript. 
BO visualized the experimental results. 
JN modeled the ASP formulation and processed the original ASP data. 
TdW and BO revised and edited the draft with support of TG, VJ and FH.
All authors revised the final manuscript.

\bibliography{SDP_QUBO_arxiv_1}
\bibliographystyle{alpha}

\appendix

\section{Technical Details in Problem Formulations}
\label{appendix:details_problem_instances}
\subsection{Affinity Matrix in ASP}
\label{appendix:details_problem_instances_ASP}
We provide the definition of the ASP in \eqref{ASP-1}.
In this section, we give additional details for the affinity measure.
We have that $\mathcal{M}$ is the set of material types and $\sigma\in \mathbb{R}^{|\mathcal{M}|\times |\mathcal{M}|}$ is a real, symmetric matrix with positive entries and zeros on the diagonal, which determines the affinity $\sigma_{mn}$ for every pair of materials $m$ and $n$ in $\mathcal{M}$.
More specifically, 
\begin{align*}
    \sigma = \frac{O_{mn}}{\delta_m + \delta_n-O_{mn}},
\end{align*}
where $O \in \mathbb{R}_{\geq 0}^{|\mathcal{M}|\times|\mathcal{M}|}$ is a symmetric matrix with zeros on the diagonal and $\delta_m \in \mathbb{Z}_{\geq 0}$. The number $O_{mn}$ is the number of times the materials $m$ and $n$ appear together in an order and $\delta_m$ is the total number of orders in which $m$ appears. Therefore  $\delta_m + \delta_n > O_{mn}$ for all $m,n\in \mathcal{M}$ -- otherwise the problem in Eq. \eqref{ASP-1} is degenerate.

\subsection{QUBO Reformulation}
\label{qubo_refomulation_appendix}
Below, we demonstrate the process of converting the IQP formulation of an ASP to its QUBO formulation. 
We exemplarily consider the smallest ASP instance \texttt{ASP30} with 30 materials and three aisles with a capacity of ten storage locations.
For the QUBO reformulation we follow the methodology developed by Glover et al. \cite{fred2019}.
\subsubsection{Unconstraining}
\label{sec:unconstraining}
Consider the IQP model of the ASP in \eqref{ASP-1}. 
In this step, we incorporate the constraints into the objective function. 
We can include the equality constraints \eqref{exactlyoneaisle} directly, while inequality constraints \eqref{aislecapcity} require slack variables to convert them into equalities before incorporation.
For each 
aisle $A_j$ with $j=1,\dots,k$ we introduce a \struc{slack variable $s_j$} in $\mathbb{Z}_{\geq 0}$. 
Moreover, let $\struc{\lambda} \in \mathbb{Z}$ be the \struc{scalar penalization factor}.
We rewrite the IQP \eqref{ASP-1} as
\begin{alignat}{3}
        &\min &&\quad \sum_{m,n\in\mathcal{M}}\sum_{\underset{i\not = j}{i,j=1}}^k\sigma_{mn}x_{mi}x_{nj} +
    \lambda \sum_{m \in \mathcal{M}}\left(\sum_{j=1}^kx_{mj} - 1\right)^2&&\label{unconstraining}\\ 
    & && \quad + \lambda \sum_{j = 1}^k\left(\sum_{m\in\mathcal{M}}x_{mj} + s_j - |A_j|\right)^2&&\notag\\
    &\textnormal{subject to } && \quad x_{mj} \in \{0,1\}, &&\text{for all } m \in \mathcal{M}, \notag\\
    &&&&&\text{and } j=1,\dots,k
    \notag\\
    & && \quad 0 \leq s_{j} \leq |A_j|,  &&\text{for all } j=1,\dots,k.\notag 
\end{alignat}

Note that it is possible to introduce an integer slack variable $s_j$ for $j=1,\dots, k$ since, in the inequality constraint (\ref{aislecapcity}), the summand on the left, as well as the aisle capacity $|A_j|$ on the right, are integers. 
The scalar penalization factor $\lambda$
aims to ensure that the original constraints are satisfied, see \cref{sec:penalization} for details.
Note that the quadratic penalty terms in the objective function in (\ref{unconstraining}) amplify larger deviations when the constraints are violated. 

Now, in the example of \texttt{ASP30}, since the number of aisles is three, the total number of integer slack variables is also three.

\subsubsection{Binarization}
\label{sec:binarization}
The above step removes the constraints \eqref{exactlyoneaisle} and \eqref{aislecapcity} from the original ASP formulation and introduces integer variables $s_j\in \mathbb{Z}_{\geq0}$, for $j=1,\dots,k$, where $k$ is the number of aisles.
In this step, we encode each of the integer slack variables with binary variables. 

Recall, that $|A_j|$ is the aisle capacity of an aisle $A_j$. 
In what follows, let $r=\lfloor\log_{2}(|A_j|)\rfloor$.
For each aisle $A_j$ with $j = 1,...,k$, we define the corresponding slack variable as
\begin{align}\label{binarize:slack}
   s_j = \sum_{s=1}^{r+1} c_{s}^{(j)}y_{s}^{(j)}, 
\end{align}
via introducing (for $s~=~1,...,r+1$) \struc{auxiliary binary variables $y_s^{(j)}$} $\in \{0,1\}$ and \struc{coefficients} $\struc{c_s^{(j)}}\in \mathbb{N}$.
Furthermore, we can define the coefficients $c_s^{(j)}$, for all $s=1, \dots, r+1$ and $j=1,\dots,k$ in detail via 
\begin{equation*}
c_{s}^{(j)} = \begin{cases}
  2^{s-1}, &\text{if } s=1,\dots, r \\
  |A_j| - (2^{r}-1), &\text{if } s = r+1.
\end{cases}\end{equation*} 

In the example of \texttt{ASP30}, we have $|A_j| = 10$ for $j=1,2,3$, since we have three aisles with ten storage locations per aisle.
This requires a total of four binary variables
for each aisle, 
leading to a total of 12 auxiliary binary variables.

With the definition in \eqref{binarize:slack}, we replace the slack variables in \eqref{unconstraining}.
This leads to
\begin{alignat}{3}
 &\min &&\quad \sum_{m,n\in\mathcal{M}}\sum_{\underset{i\not = j}{i,j=1}}^k\sigma_{mn}x_{mi}x_{nj} + 
    \lambda \sum_{m\in\mathcal{M}}\left(\sum_{j=1}^kx_{mj} - 1\right)^2&&\label{binarizing2} \\
    &&&\quad + 
    \lambda \sum_{j=1}^k\left(\sum_{m\in\mathcal{M}}x_{mj} + \sum_{s=1}^{r+1} c_{s}^{(j)}y_{s}^{(j)} - |A_j|\right)^2&&\notag\\
     &\textnormal{subject to } && \quad x_{mj} \in \{0,1\}, && \text{ for all } m \in \mathcal{M},\notag\\
     &&&&&\text{ and }j=1,\dots,k\notag\\
     & &&\quad  y_{s}^{(j)} \in \{0,1\}, \qquad \qquad && \text{ for all } s =1,\dots,r+1,\notag\\
     &&&&&\text{ and } j=1,\dots,k \notag.
\end{alignat}

Finally, we reach the QUBO formulation in \eqref{eq:QUBO}, 
via expanding quadratic penalty terms and omitting the constant term, which does not affect optimization.

Considering the example \texttt{ASP30}, we reformulate the original problem \eqref{ASP-1} as a QUBO with 96 binary variables.
These consists of 12 auxiliary binary variables and 84 for the assignment of each material $m\in \mathcal{M}$ to an aisle $A_j$ with $j=1,\dots,k$. 
Since two materials $m,n\in \mathcal{M}$ have affinity $\sigma_{mn} = 0$, the binary variables reduce from 90 to 84. 

Note that, as we demonstrate in this subsection,  
the number of binary variables required to reformulate an IQP, and similarly an IP, into a QUBO grows logarithmically with the upper bound of the slack variables, which in this case is the aisle capacity of each aisle.

\subsubsection{Penalization}
\label{sec:penalization}
It is nontrivial to select the appropriate penalization factor $\lambda\in\mathbb{Z}$ in \eqref{unconstraining}, respectively \eqref{binarizing2}, for the reformulation. 
However, knowledge about the range of the objective function value can be very useful in this regard. 
As explained in \cite{fred2019} and \cite{10.1007/s10878-014-9734-0}, there is a ``Goldilocks region" where penalty values strike a balance and work effectively. 
This region is typically broad, offering some flexibility. 
To estimate an appropriate penalization factor, a preliminary analysis of the original model can help provide a rough approximation of the objective function's value. 
A good starting point is to set the penalty between 75\% and 150\% of this estimate, then iteratively adjust based on feasibility until an acceptable solution is found.

\section{Further Certificates of Nonnegativity beyond SOS-SDP}
\label{sec:nonnegativity_certificates}
As pointed out, QUBOs fall into the broader range of polynomial optimization problems (on the boolean hypercube).
There exist various other certificates of nonnegativity beyond SOS that are tractable via convex optimization problems.
These include for example \struc{sums of nonnegative circuit polynomials (SONC)} certificates, which are also referred to \struc{sums of AM/GM exponentials (SAGE)}, which can be computed via relative entropy programs; see e.g., \cite{Iliman:deWolff:Circuits,Chandrasekaran:Shah:SAGE}.
On the boolean hypercube, another prominent hierarchy is the Sherali-Adams hierarchy, which leads to a linear programming formulation \cite{Sherali:Adams}.
Furthermore, recently a new certificate named \struc{sums of copositive fewnomials} was suggested \cite{Averkov:Scheiderer}.
We decided against including these certificates in the experimental comparison carried out in this work first due to the complexity of the endeavor, second since at least in some of the cases proper implementations are not available, and third since we aim to focus on a comparison of SDP-based methods with Hamiltonian Updates.

Finally, a very natural certificate to consider here would be \struc{scaled diagonally dominant sum of squares (SDSOS)}, which are tractable by SOCPs \cite{Ahmadi:Majumdar}. SDSOS are hence a special case of SOS, which are, however, easier to compute (they happen moreover to be SONCs).
Unfortunately, as it was confirmed to us by one of the developers, no implementation suitable for the experiments carried out in this paper is available at the moment.

\section{Converting QUBOs to Ising formulation}
\label{subsection:binary_conversion}

We convert our initial QUBO problem $\min_{\Vector{b}} \Vector{b}^TQ\Vector{b}$, with $\Vector{b}\in\{0,1\}^n$ and $Q\in \mathbb{R}^{n\times n}$, see \eqref{eq:QUBO}, to $\min_{\Vector{x}} \Vector{x}^TC\Vector{x}$ with $\Vector{x}\in\{-1,1\}^{n+1}$ and $C\in\mathbb{R}^{(n+1) \times (n+1)}$. This formulation is also known as the \struc{Ising formulation} used for Hamiltonian Updates, see \eqref{eq:ising}, as well as an alternative input for the polynomial optimization software.
In order to construct $C$, we insert an additional variable $x_0$ that we constrain to be $x_0 = 1$.
We convert between the vectors $\Vector{b}$ and $\Vector{x}$ via $b_i = (x_i+1)/2$ for all $i=1,\dots, n$.
The new polynomial $f(1,\Vector{x})$ given by the new cost matrix $C\in \mathbb{R}^{(n+1)\times (n+1)}$ looks as follows. 
Recall that $Q$ is a symmetric matrix.
Thus, we have
\begin{align*}
    f(1,\Vector{x}) &= [1,\Vector{x}]^TC[1,\Vector{x}]\\ &= \begin{bmatrix}
        1&x_1&x_2\hdots&x_n
    \end{bmatrix} \frac{1}{4} 
    \begin{bmatrix}
    \begin{array}{c|ccc}
 \sum_{i,j=1}^n q_{ij} & \sum_{i=1}^nq_{i1} & \hdots & \sum_{i=1}^nq_{in}\\
  \hline
  \sum_{i=1}^nq_{i1} & q_{11} & \hdots & q_{1n} \\
        \vdots & \vdots &\ddots & \vdots \\
        \sum_{i=1}^nq_{in} & q_{1n} & \hdots & q_{nn}
    \end{array}
    \end{bmatrix}
    \begin{bmatrix}
        1\\x_1\\x_2\\ \vdots \\ x_n
    \end{bmatrix} \\
    &= \frac{1}{4}\Big(\Vector{x}^TQ\Vector{x} + 2\sum_{j=1}^n x_j\sum_{i=1}^n q_{ij} + \sum_{i,j}^n q_{ij} \Big).
\end{align*}
The advantage of using $\{-1,1\}$ as binary representation instead of $\{0,1\}$ is that all monomials of the form $q_{ii}x_i^2$ become directly constant, i.e., $q_{ii}x_i^2 = q_{ii}$ for $x_i = \pm 1$.
This directly reduces the number of monomials, in particular all $n$ monomials corresponding to the diagonal of $C$ are now constant.

\section{Technical Details in SDP Relaxations}
\label{app:tech_details}
\subsection{Time Limits}
In \cref{tab:timelimits} we display the time limits we enforce for the SDP relaxations.
For a single parameter choice, we give a time limit of the number of instances in each set times one hour for ASP instances, and times two hours for OVRP instances.
For example, since there are 14 QUBOs corresponding to \texttt{ASP30}, we enforce a time limit of 14 hours for one run with a specific parameter set through all \texttt{ASP30} instances.
Exemplarily, we rerun and provide outcomes of computations that were exceeding the initial time limit.
	\begin{table}[h!]
		\begin{center}
			\begin{tabular}{|l||c|c|c||c|}
            \hline
				Instance Set &	 \texttt{ASP30} &  \texttt{ASP60} &  \texttt{ASP90} & OVRP  \\
                \hline
               Time Limit in hours &14 &4 &3 &4 \\
				 \hline
			\end{tabular}
		\end{center}
        \caption{Imposed \textbf{time limits} in hours for one run through the entire set of instances with each parameter choice.}
        \label{tab:timelimits}
	\end{table}
    
\subsection{SOS Solver Settings} For the SDP solvers inside \textsc{TSSOS} \cite{TSSOSDok}, we test different tolerances in the solver settings in order to achieve solutions with higher/lower precision:
\begin{itemize}
    \item For \textsc{Mosek} \cite{mosek}, we set \textit{primal feasibility tolerance}, \textit{dual feasibility tolerance}, as well as \textit{primal-dual gap tolerance} to $\{10^{-8}, 10^{-4}\}$ and refer to default parameters as set up in the documentation of \cite{TSSOS}.
    \item For \textsc{COSMO} \cite{cosmo}, we set \textit{absolute residual tolerance} and \textit{relative residual tolerance} to $\{10^{-5}, 10^{-3}\}$ and refer to default parameters as set up in the documentation of \cite{TSSOS}.
\end{itemize}
Using the \textsc{ManiSDP} solver directly for the moment version of \eqref{TUBS:lasserre}, we use the default parameters for SDPs with diagonal and unit-diagonal constraints respectively, as given in the documentation of \cite{wang2024solvinglowranksemidefiniteprograms}. 
For \textsc{SOS.jl}, we use the default settings of \textsc{Mosek} \cite{mosek} and \textsc{Loraine} \cite{loraine2023}. We also use the default settings for the solver \textsc{SCS} \cite{scs}.

\section{Additional Data}
\label{appendix:tables}
In the additional tables we provide information about the absolute differences $\DeltaQUBOabs$, $\DeltaSOSabs$, $\DeltaHUabs$ to the optimal values, as well as $\DeltaHUabsupper, \DeltaHUupper$ for the upper bounds achieved via the randomized rounding procedure in HU. 
For reference, the optimal values are $\Zopt~\approx~246.78$ for \texttt{ASP30}, $\Zopt\approx820.08$ for \texttt{ASP60} and $\Zopt\approx1373.15$ for \texttt{ASP90}, $\Zopt \approx 6254.19$ for   \texttt{OVRP11} and $\Zopt\approx 3036.54$ for \texttt{OVRP15}, see~\cref{tab:MIP_results}.

\subsection{SOS-SDP} 
In addition to \cref{tab:ASP_running_times}, we give the \textbf{average} of the \textbf{absolute difference} $\DeltaSOSabs$ of the lower bound $\zSOS$ to the optimal value using SOS-SDP methods for the ASP instances in \cref{tab:ASP_running_times_abs}.

	\begin{table}[h!]
		\begin{center}
			\begin{tabular}{|l|c|c|c|c|c|c|}\hline
				&	\multicolumn{2}{c|}{\texttt{ASP30}} & \multicolumn{2}{c|}{\texttt{ASP60}}  & \multicolumn{2}{c|}{\texttt{ASP90}} \\
				&\multirow{2}*{} \multirow{2}*{Average} & \multirow{2}*{Average} & \multirow{2}*{Average} & \multirow{2}*{Average}& \multirow{2}*{Average}& \multirow{2}*{Average}  \\[+4pt] 
                & $\DeltaSOSabs$ & $\tSOS$ &$\DeltaSOSabs$ &$\tSOS$ &$\DeltaSOSabs$ & $\tSOS$ \\
                \hline
                \textsc{ManiSDP} & $491.1$& $1.45 s$& $2468.6$ &  $0.82 s$& $3308.6$ & $1.45 s$\\
				\hline 
                \textsc{SOS.jl} w. \textsc{Mosek}  & $491.1$& $9.9 s$& $2468.6 $ &  $221.25 s$& $3308.5$ & $8256.9 s$\\
				\hline   
                \textsc{SOS.jl} w. Loraine  & $491.1$& $670.66 s$& $\dagger$ &  $\dagger$& $\dagger$ & $\dagger$\\
				\hline 
                \textsc{TSSOS} w. \textsc{Mosek} & $491.1$& $3.24 s$& $2468.6$ &  $112.42 s$& $\dagger$ & $\dagger$\\
				\hline 
                \rowcolor{gray!20}\textsc{TSSOS} w. \textsc{Mosek} & $32.72$& $28.46 s$& $174.3$ &  $912.47s$& $\dagger$ & $\dagger$\\ \hline 
			\end{tabular}
		\end{center}
        \caption{\textbf{Results of SOS-SDP relaxations} of QUBOs referring to \texttt{ASP30}, \texttt{ASP60} and \texttt{ASP90}. We display the \textbf{average} of the \textbf{absolute difference} $\DeltaSOSabs$ of the lower bound $\zSOS$ to the optimal value. The gray row refers to the second order Lasserre relaxation, the other rows refer to the first order relaxation. A dagger '$\dagger$' indicates that no result is available due to solver break down or exceeding maximal running time. For \textsc{Mosek}, we display the result with the default tolerance tol=$10^{-8}$.  For \textsc{SOS.jl} with \textsc{Mosek} the total solving time of all \texttt{ASP90} instances exceeds our set time limit. We display a rerun for these instances without time limit. }
        \label{tab:ASP_running_times_abs}
	\end{table}

\subsection{HU}
In addition to \cref{tab:VRPresults_GPU,tab:ASP60_90results_GPU_example}, we show the \textbf{absolute difference} of the \textbf{lower bound} $\DeltaHUabs$ to $\Zopt$, and the absolute and relative differences of the \textbf{upper bound} $\DeltaHUabsupper$ and $\DeltaHUupper$ for OVRP, as well as selected ASP instances, using HU in \cref{tab:VRPresults_GPU_upperbounds,tab:ASP60_90results_GPU_example_upper}, respectively. 

	\begin{table}[h!]
		\begin{center}
			\begin{tabular}{|l|c|c|c|c|c|c|}\hline
				&	\multicolumn{3}{c|}{\texttt{OVRP11}} & \multicolumn{3}{c|}{\texttt{OVRP15}} \\
				
				& $\DeltaHUabs$ & $\DeltaHUabsupper$& $\DeltaHUupper$& $\DeltaHUabs$ & $\DeltaHUabsupper$& $\DeltaHUupper$ \\
				\hline
				HU$^*$\tiny{$\epsilon = \num{e-2}$} &$ \num{5.0e10}$ & $\num{1.4e9}$ & $ \num{2.1e5}$ & $ \num{1.7e12}$ & $\num{2.5e9}$
					& $ \num{8.4e5}$ \\
				HU$^*$\tiny{$\epsilon = \num{e-3}$} &$ \num{4.7e9}$ & $\num{1.3e9}$ & $\num{2.1e5}$ & $\num{1.9e10}$ & $\num{2.6e9}$ & $\num{8.5e5}$\\
                HU$^*$\tiny{$\epsilon = \num{e-4}$} &$ \num{4.4e8}$ & $\num{3.7e7}$ & $\num{5.9e3}$ & $\num{1.7e9}$ & $\num{6.7e7}$ & $\num{2.2e4}$\\
                \hline
			\end{tabular}
		\end{center}
        \caption{\textbf{Results of SDP relaxations via HU} of QUBOs referring to \texttt{OVRP11} and \texttt{OVRP15} on \textbf{GPU}. We display the absolute difference of the \textbf{lower bound} $\DeltaHUabs$ to $\Zopt$, and the absolute and relative differences of the \textbf{upper bound} $\DeltaHUabsupper$ and $\DeltaHUupper$. The '$*$' indicates that the results are obtained via GPU.}
        \label{tab:VRPresults_GPU_upperbounds}
	\end{table}
    
	\begin{table}[h!]
		\begin{center}
			\begin{tabular}{|l|c|c|c|c|c|c|}\hline
				&	\multicolumn{3}{c|}{\texttt{ASP60-850}} &	\multicolumn{3}{c|}{\texttt{ASP90-1400}} \\
				& $\DeltaHUabs$ & $\DeltaHUabsupper$ & $\DeltaHUupper$ &  $\DeltaHUabs$ & $\DeltaHUabsupper$ & $\DeltaHUupper$  \\
                \hline
                HU$^*$\tiny{$\epsilon = \num{e-2}$} & $\num{4.4e6}$ &$\num{1.4e5}$ & $172.3$& $\num{2.7e7}$ &$\num{5.3e5}$ & $\num{386.1}$ \\
                 HU$^*$\tiny{$\epsilon = \num{e-3}$} &$\num{6.6e5}$ &$\num{2.9e4}$ &$35.9$ & $\num{3.7e6}$ &$\num{7.7e4}$ & $55.8$\\
                   HU$^*$\tiny{$\epsilon = \num{e-4}$} &$\num{6.1e4}$ &$\num{1.3e4}$ &$15.9$ & $\num{3.5e5}$ &$\num{3.6e4}$ & $25.9$\\
                    \hline
			\end{tabular}
		\end{center}
        \caption{\textbf{Results of SDP relaxations via HU} of QUBOs referring to \texttt{ASP60-850} and \texttt{ASP90-1400}. We display the absolute difference of the \textbf{lower bound} $\DeltaHUabs$ to $\Zopt$, and the absolute and relative differences of the \textbf{upper bound} $\DeltaHUabsupper$ and $\DeltaHUupper$. The '$*$' indicates that the results are obtained via GPU.}
        \label{tab:ASP60_90results_GPU_example_upper}
	\end{table}
    
\subsection{QUBO Solver Results}
In addition to \cref{tab:QUBO_example_optimality}, in \cref{tab:QUBO_example_optimality_abs_difference} we display the \textbf{absolute difference} $\DeltaQUBOabs$ to the optimal value $\Zopt$, solving the QUBO formulation with \textsc{GUROBI}.
\label{sec:appendix_qubo_solver_results}
\begin{table}[H]
\begin{center}
    \begin{tabular}{|l||c|c||c|c||c|c|}
    \hline
        Instance & Time Limit & $\DeltaQUBOabs$ & Time Limit & $\DeltaQUBOabs$ & Time Limit & $\DeltaQUBOabs$ \\
        \hline
         \texttt{ASP30}-$950$ & $10s$ & $0$ &$120s$ &$0$ &$600s$ &$0$\\
         \texttt{ASP60}-$900$ & $10s$ & $240.3$ &$120s$ & $240.3$&$600s$ &$240.3$\\
         \texttt{ASP90}-$1450$ & $10s$ & $1020.6$ &$120s$ &$270.9$ &$600s$ &$262.9$\\
        \hline
    \end{tabular}
\end{center}
\caption{\textbf{Absolute difference} $\DeltaQUBOabs$ of the \textbf{upper bound} $\zQUBO$ to the optimal value $\Zopt$ for selected \textbf{QUBO}s using \textsc{GUROBI} with different time limits. For \texttt{ASP30}-950 we find the optimal value. The provided time limits are upper bounds on the running time to reach the respective objective value $\zQUBO$. The computation is stopped at the provided time limit.}
\label{tab:QUBO_example_optimality_abs_difference}
\end{table}

\subsection{Additional Quantum Annealer Results}
\label{sec:dwave_results}
In addition to \cref{sec:mip_solutions} and \cref{tab:QUBO_example_optimality_abs_difference}, where we solve the QUBO formulations via the classic and commercial solvers \textsc{CPLEX} and \textsc{GUROBI}, we exemplarily give an overview on solving them via the quantum annealer from D-Wave \cite{DWaveLeap}, using the \textsc{QUBO\_SOLVE} tool from \textsc{GAMS}, see \cref{sec:ex_setup}.
We emphasize, that quantum annealing is not the main focus of our study, we provide these results for reference.
Moreover, we point out, that these results are not easily reproducible as the free developer access is no longer available.
We see in \cref{tab:dwave_results} that, similar to solving via the classical solvers, we find feasible solutions for the ASP instances, but the solutions for the OVRP instances are infeasible.
\begin{table}[h!]
		\begin{center}
    \begin{tabular}{|l|c|c|c|c|c|}
    \hline
        Instance & \texttt{ASP30}-$450$ & \texttt{ASP60}-$900$ & \texttt{ASP90}-$1450$ & \texttt{OVRP11} & \texttt{OVRP15} \\\hline\hline
        $\DeltaQUBOabs$ &$0$ &$65.4$ & $843.4$& infeasible & infeasible\\
        $\DeltaQUBO$ &$0$ &$0.08$ & $0.61$& infeasible & infeasible\\
        Time Limit & $120s$& $120s$& $120s$& $60s$& $60s$\\\hline
    \end{tabular}
\end{center}
\caption{\textbf{Absolute difference} $\DeltaQUBOabs$ and \textbf{relative difference} $\DeltaQUBO$ of the \textbf{upper bound} $\zQUBO$ to the optimal value $\Zopt$ for selected \textbf{QUBO}s using \textbf{quantum annealing} by D-Wave with indicated time limits. For \texttt{ASP30}-450 we find the optimal value. We stop the computation at the given time limit. Note that we also find the optimal value for \texttt{ASP30}-450 for a time limit of 30s and 60s.}
        \label{tab:dwave_results}
	\end{table}

\end{document}